\documentclass[final,12pt]{msml2021} 

\def\mathrlap#1{\text{\hbox to 10pt{$\mathsurround=0pt#1$\hss\hss}}}
\usepackage{latexsym}
\usepackage{indentfirst}
\usepackage{placeins}
\usepackage{booktabs}
\usepackage{algorithm}
\usepackage{algorithmic}
\usepackage{multirow}
\usepackage{color}
\usepackage{diagbox}

\newtheorem{prop}[theorem]{Proposition}

\newtheorem{defn}[theorem]{Definition}

\newcommand{\grad}{\nabla}

\newcommand{\x}{\mathbf{x}}
\newcommand{\E}{\mathbb{E}}

\title[Semigroup method for committor functions]{A semigroup method for  high dimensional
  committor functions based on neural network}

\usepackage{times}

\msmlauthor{%
 \Name{Haoya Li} \Email{lihaoya@stanford.edu}\\
 \addr Department of Mathematics, Stanford University, Stanford, CA 94305, USA
 \AND
 \Name{Yuehaw Khoo} \Email{ykhoo@statistics.uchicago.edu}\\
 \addr Department of Statistics and the College, The University of Chicago, Chicago, IL 60637, USA
 \AND
 \Name{Yinuo Ren} \Email{renyinuo@pku.edu.cn}\\
 \addr School of Mathematical Sciences, Peking University, Beijing, China 
 \AND
 \Name{Lexing Ying} \Email{lexing@stanford.edu}\\
 \addr Department of Mathematics and ICME, Stanford University, Stanford, CA
	94305, USA
}

\makeatletter
 \let\Ginclude@graphics\@org@Ginclude@graphics
\makeatother
\begin{document}

\maketitle

\begin{abstract}%
  This paper proposes a new method based on neural networks for computing the high-dimensional
  committor functions that satisfy Fokker-Planck equations. Instead of working with partial
  differential equations, the new method works with an integral formulation based on the semigroup
  of the differential operator. The variational form of the new formulation is then solved by
  parameterizing the committor function as a neural network. There are two major benefits of
    this new approach. First, stochastic gradient descent type algorithms can be applied in the
    training of the committor function without the need of computing any mixed second order
    derivatives. Moreover, unlike the previous methods that enforce the boundary conditions through
    penalty terms, the new method takes into account the boundary conditions automatically.
  Numerical results are provided to demonstrate the performance of the proposed method.

\end{abstract}


\section{Introduction}\label{sec:intro}
Understanding rare transition events between two states is important for studying the behavior of
stochastic systems in physics, chemistry, and biology. One important method to describe the
transition events is the transition path theory
\citep{vanden2010transition,weinan2006towards,lu2015reactive}, and the central object in the
transition path theory is the committor function. Assume that the transition between two states is
governed by the overdamped Langevin process
\begin{equation}\label{langevin} 
  d \x_{t}=-\nabla V\left(\x_{t}\right) d t+\sqrt{2 \beta^{-1}} d \mathbf{w}_{t}, 
\end{equation}
where $\x_t\in \Omega \subset \mathbb{R}^d$ is the state of the system,
$V:\mathbb{R}^d\rightarrow \mathbb{R}$ is a smooth potential function, $\beta = 1/T$ is the inverse
of the temperature $T$, and $\mathbf{w}_t$ is the standard $d$-dimensional Brownian motion.  For two
given simply connected domains $A$ and $B$ in $\Omega$ with smooth boundaries, the committor function $q(\mathbf{x})$ is defined as
\begin{equation}
  q(\x) = \mathbb{P}(\tau_B< \tau_A\mid \x_0 = \x), 
\end{equation}
where $\tau_A$ and $\tau_B$ are the hitting times for the sets $A$ and $B$, respectively. Many
statistical properties of the transition are encoded in the committor function.

The committor function satisfies the Fokker-Planck equation (also known as the steady-state backward
Kolmogorov equation)
\begin{equation}\label{committor}
  (-1/\beta \Delta + \nabla V\cdot \grad) q = 0 \text { in } \Omega \backslash (A \cup B), \quad \left.q(\x)\right|_{\partial A}=0, \quad \left.q(\x)\right|_{\partial B}=1. 
\end{equation}
If the potential $V$ is confining, i.e. $V(\x)\rightarrow \infty$ as
$|\x|\rightarrow\infty$ and the partition function $Z_\beta = \int_{\mathbb{R}^d} \exp(-\beta
V(\x))d \x<\infty$ for any $\beta>0$, then
\begin{equation}\label{rhobeta}
  \rho(\x) = \frac{1}{Z_\beta} \exp(-\beta V(\x))
\end{equation}
is the equilibrium distribution of the Langevin dynamics \eqref{langevin}.  Under this condition,
the Langevin process \eqref{langevin} is ergodic \citep{pavliotis2014stochastic}, which enables the use
of Monte Carlo methods to sample from the equilibrium distribution. For the rest of the paper we
will always assume that $V$ is confining and $\rho$ exists.

One major difficulty in solving \eqref{committor} is the curse of dimensionality. Various methods
have been proposed under the assumption that the transition from $A$ to $B$ is concentrated in a
quasi-one dimensional tube or low dimensional manifold. For example, the finite temperature string
method \citep{weinan2005finite,vanden2009revisiting} approximates isosurfaces of the committor
function with hyperplanes normal to the most probable transition paths, and updates the transition
paths together with the isocommittor-surfaces. Diffusion map \citep{coifman2008diffusion} solves for
$q$ on a set of points by applying point cloud discretization to the operator
$L=(-1/\beta\Delta+\nabla V \cdot \grad)$.  In order to obtain better convergence order,
\citep{lai2018point} improves on diffusion map by discretizing $L$ using a finite element method on
local tangent planes of the point cloud.

More recently, the method proposed in \citep{khoo2019solving} works with the variational form of the
Fokker-Planck equation
\begin{equation}\label{ritz}
  \underset{q}{\operatorname{argmin}} \int_{\Omega \backslash (A \cup B)}|\nabla
  q(\x)|^{2} \rho(\x) d \x,\left.\quad
  q(\x)\right|_{\partial A}=0,\left.\quad q(\x)\right|_{\partial B}=1, 
\end{equation}
and parameterizes the high-dimensional committor function $q(\cdot)$ by a neural network (NN)
$q_{\theta}(\cdot)$. The main advantage of working with an optimization formulation is that
under rather mild conditions stochastic gradient descent (SGD) type algorithms can avoid the
  saddle points and converge efficiently to at least local minimums, as shown for example in
  \cite{jin2019nonconvex}.  The boundary conditions can be enforced by additional penalty terms as
in
\begin{equation}\label{weakpenalty}
\begin{aligned}
  \underset{\theta \in \mathbb{R}^{p}}{\operatorname{argmin}} \int_{\mathrlap{\Omega \backslash (A \cup B)}}\left|\nabla q_{\theta}(\x)\right|^{2} \rho(\x) d \x
  +\tilde{c} \int_{\partial A} q_{\theta}(\x)^{2} d m_{\partial A}(\x)+\tilde{c} \int_{\partial B}\left(q_{\theta}(\x)-1\right)^{2} d m_{\partial B}(\x), 
\end{aligned}
\end{equation}
where $m_A$ and $m_B$ are measures supported on $\partial A$ and $\partial B$.  The integral in
\eqref{weakpenalty} is approximated by a Monte-Carlo method where the samples are generated
according to the stochastic process \eqref{langevin}. Satisfactory numerical results are obtained
with this method even in high dimensional cases where traditional methods like finite-element-type
methods are intractable. However, because the loss depends on $|\nabla q_{\theta}|^2$, obtaining the
gradient of the loss with respect to $\theta$ requires an inconvenient second order derivative
computation.

The contribution of this paper is two-fold. First, we propose a new method that removes
  the dependence on the second order derivatives by working with an integral formulation based on
  the semigroup of \eqref{langevin}. Second, the boundary conditions are treated automatically by
  the semigroup formulation rather than relying solely on the penalty terms, which makes it easier
  to tune the penalty coefficients. The paper is organized as follows. Section~\ref{sec:alg}
describes the new formulation and the neural network approximation. Section~\ref{sec:lazy} analyzes
this new method under the so-called lazy-training regime. Finally, the performance of our method is
examined numerically in Section~\ref{sec:num}.

\section{Proposed method}\label{sec:alg}

\subsection{A new variational formulation}\label{sec:consist}

Consider the Langevin process starting from a point $\x\in\Omega$
\begin{equation}\label{eq:trajectory}
  \begin{aligned}
    d \x_{t} &= -\nabla V\left(\x_{t}\right) d t+\sqrt{2 \beta^{-1}} d \mathbf{w}_{t}, \\
    \x_0 &= \x.  
  \end{aligned}
\end{equation}
Let $\tau_A$ and $\tau_B$ be the stopping time of the process hitting $\partial A$ and $\partial B$,
respectively. Similarly, define $\tau \equiv \tau_{A\cup B} = \min(\tau_A,\tau_B)$ be the hitting
time of $A\cup B$. For a fixed small time step $\delta>0$, we define the operator $P$ as follows:
\begin{equation}\label{eq:beforesplit}
  (Pf)(\x):=\E^{\x}\left(f\left(\x_{\tau\wedge \delta}\right)\right),
\end{equation}
where $\E^{\x}$ is the expectation taken with respect to the law of the process
\eqref{eq:trajectory}.

\begin{prop}\label{prop:semi}
  When $\nabla V$ is bounded and Lipschitz continuous on $\mathbb{R}^d$, the solution to the
  committor function \eqref{committor} satisfies the following semigroup formulation
  \begin{equation}\label{eq:justifiedreplace}
    q(\x) = (Pq)(\x) \quad\text{in}\quad \Omega \backslash (A \cup B), \quad
    q|_{\partial A} = 0, \quad q|_{\partial B} = 1.
  \end{equation}
\end{prop}
The proof of Proposition \ref{prop:semi} is provided in Appendix~\ref{ap:semi}. The main advantage
of working with \eqref{eq:justifiedreplace} over \eqref{committor} is that it does not contain any
differential operator.

For notational convenience, we introduce a function $r:\partial A\cup\partial B\rightarrow\mathbb{R}$
with $r(\x)|_{\partial A}=0$ and $r(\x)|_{\partial B}=1$. With this definition, the
boundary condition of \eqref{eq:justifiedreplace} is simply $q=r$ on $\partial A\cup\partial B$. In order
to introduce the variational formulation, $Pq$ can be split into two parts as follows:
\begin{equation}\label{eq:split}
  \begin{aligned}
    (P q)(\x) & = \E^{\x}\left(q\left(\x_{\tau\wedge \delta}\right)\right) = \E^{\x}\left(q\left(\x_{\tau\wedge \delta}\right)\mathbf{1}_{\{\delta<\tau\}}\right) + \E^{\x}\left(q\left(\x_{\tau\wedge \delta}\right)\mathbf{1}_{\{\delta\geq \tau\}}\right)\\
    & = \E^{\x}\left(q\left(\x_{\delta}\right)\mathbf{1}_{\{\delta<\tau\}}\right) + \E^{\x}\left(q\left(\x_{\tau}\right)\mathbf{1}_{\{\delta\geq \tau\}}\right)\\
    &= \E^{\x}\left(q\left(\x_{\delta}\right)\mathbf{1}_{\{\delta<\tau\}}\right) + \E^{\x}\left(r\left(\x_{\tau}\right)\mathbf{1}_{\{\delta\geq \tau\}}\right), 
  \end{aligned}
\end{equation}
where the last equality results from the fact that $\x_{\tau}\in \partial A\cup\partial B$ and $q=r$ on $\partial A\cup\partial B$. 

We denote the first part of \eqref{eq:split} as
\begin{equation}\label{eq:pi}
  (P^i q)(\x) \equiv
  \E^{\x}\left(q(\x_{\tau\wedge\delta})\mathbf{1}_{\{\delta<\tau\}}\right) =
  \E^{\x}\left(q(\x_{\delta})         \mathbf{1}_{\{\delta<\tau\}}\right),
\end{equation}
where the superscript $i$ stands for the interior contribution and the second part of \eqref{eq:split} as
\begin{equation}
  (P^b r)(\x) \equiv
  \E^{\x}\left(r(\x_{\tau\wedge\delta})\mathbf{1}_{\{\delta\geq\tau\}}\right)
  =
  \E^{\x}\left(r(\x_{\tau})          \mathbf{1}_{\{\delta\geq\tau\}}\right),
\end{equation}
where the superscript $b$ stands for the boundary contribution. With these definitions, we can rewrite
\eqref{eq:justifiedreplace} compactly as
\begin{equation}\label{eq:nobd}
  (I-P^i)q(\x)-(P^b r)(\x) = 0, 
\end{equation}
where $I$ is the identity operator. The following result states that $P^{i}$ is symmetric on the
Hilbert space $L_{\rho}^2(\Omega\backslash (A\cup B))$.

\begin{prop}\label{prop:sym}
$P^i$ is a symmetric operator on $L_{\rho}^2(\Omega\backslash (A\cup B))$, in other words, $\langle u,
  P^i v\rangle_{\rho} = \langle P^i u, v\rangle_{\rho}$, where $\langle \cdot, \cdot\rangle_{\rho}$
  denotes the inner product of the Hilbert space $L_{\rho}^2(\Omega\backslash (A\cup B))$.
\end{prop}
For the proof of this proposition, see Appendix~\ref{ap:sym}. 

Based on Proposition~\ref{prop:sym}, we are now ready to propose the following variational
formulation for \eqref{eq:nobd}

\begin{equation}\label{eq:variational}
  \begin{aligned}
    \min_q \frac{1}{2}\int_{\mathrlap{\Omega \backslash (A \cup B)}} q(\x)\left((I-P^i)q(\x)\right)
    \rho(\x) d \x - \int_{\mathrlap{\Omega \backslash (A \cup B)}} q(\x)P^b r(\x)\rho(\x) d \x
  \end{aligned}
\end{equation}
The two formulations \eqref{eq:variational} and \eqref{eq:nobd} share the same solution as shown
below. Let $q^*$ be the solution to the variational problem \eqref{eq:variational} and $q(\x,
\epsilon) = q^*(\x)+\epsilon\eta(\x)$, where $\eta$ is continuous with compact support. By taking derivative with respect to $\epsilon$ at $\epsilon = 0$ we obtain
\begin{equation}\label{eq:varderiv}
  \begin{aligned}
    0 &= \frac{\partial}{\partial \epsilon} \left.\left( \frac{1}{2}\int_{\mathrlap{\Omega \backslash (A \cup
        B)}}q(\x, \epsilon)\left((I-P^i)q(\x, \epsilon)\right) \rho(\x) d \x - \int_{\mathrlap{\Omega \backslash (A \cup
        B)}} q(\x, \epsilon)P^b r(\x)\rho(\x) d \x\right)\right|_{\epsilon = 0}\\
    &=\int_{\mathrlap{\Omega \backslash (A \cup B)}} \eta(\x)\left((I-P^i)q^*(\x)\right)\rho(\x)d\x
    -\int_{\mathrlap{\Omega \backslash (A \cup B)}} \eta(\x)P^b r(\x)\rho(\x) d
    \x\\
    &=\int_{\mathrlap{\Omega \backslash (A \cup B)}}
    \eta(\x)\left((I-P^i)q^*(\x)-P^b r(\x)\right)\rho(\x) d \x,
  \end{aligned}
\end{equation}
where the second equality uses Proposition~\ref{prop:sym}. Since this is true for any continuous
$\eta$ with compact support and $\rho(\x)>0$, we conclude that $(I-P^i)q^*(\x)-P^b r(\x) = 0$ on
$\Omega \backslash (A \cup B)$.

Plugging in the definitions of $P^i$ and $P^b$ into \eqref{eq:variational} leads to
\begin{equation}\label{eq:removegrad}
  \begin{aligned}
    \min_q ~\frac{1}{2}\int_{\mathrlap{\Omega \backslash (A \cup
        B)}}q(\x)\left(q(\x)-\E^{\x}\left(q\left(\x_{\delta}\right)\mathbf{1}_{\{\delta<\tau\}}\right)\right)
    \rho(\x) d \x - \int_{\mathrlap{\Omega \backslash (A \cup B)}}
    q(\x)\E^{\x}\left(r\left(\x_{\tau}\right)\mathbf{1}_{\{\delta\geq\tau\}}\right)\rho(\x)
    d \x.
  \end{aligned}
\end{equation}
It is clear from this formulation that there is no need for taking gradient of $q(\x)$ with respect to $\x$.

Although there is no need to enforce the boundary conditions $q|_{\partial A} = 0$ and $q|_{\partial
  B} = 1$ explicitly in \eqref{eq:variational}, one can still include the penalty terms that can
sometimes give a better performance
\begin{equation}\label{eq:adddpen}
  \begin{aligned}
    \min_q \frac{1}{2} & \int_{\mathrlap{\Omega \backslash (A \cup B)}} q(\x)\left((I-P^i)q(\x)\right)
    \rho(\x) d \x - \int_{\mathrlap{\Omega \backslash (A \cup B)}} q(\x)P^b r(\x)\rho(\x) d \x\\
    &+   \frac{c}{2} \int_{\partial A} q(\x)^2 dm_A(\x) + \frac{c}{2}    \int_{\partial B} (q(\x)-1)^2 dm_B(\x),
  \end{aligned}
\end{equation}
where $m_A$ and $m_B$ are measures supported on $\partial A$ and $\partial B$,  respectively, and $c>0$ is a penalty
constant.

\subsection{Nonlinear parameterization}\label{sec:nonlinear}

In order to deal with the high dimensionality of the committor function $q(\x)$, we propose to
approximate $q(\x)$ with an NN $q_\theta(\x)$, where $\theta$ stands for the NN parameters. In terms
of $\theta$, the optimization problem takes the form
\begin{equation}\label{eq:no gradient var}
  \begin{aligned}
    \min_\theta \frac{1}{2} &
    \int_{\mathrlap{\Omega \backslash (A \cup B)}} q_\theta(\x)
    \left((I-P^i)q_\theta (\x)\right) \rho(\x) d \x -
    \int_{\mathrlap{\Omega \backslash (A \cup B)}} q_\theta(\x) P^b r(\x)\rho(\x) d \x \\
    &+ \frac{c}{2} \int q_\theta(\x)^2 dm_A(\x) + \frac{c}{2} \int (q_\theta(\x)-1)^2 dm_B(\x).
  \end{aligned}
\end{equation}
When applying an stochastic gradient descent (SGD) to optimize \eqref{eq:no gradient var}, one needs to compute the derivative for
each of the terms.

\paragraph{Derivative of the first two terms.}
By the symmetric property stated in Proposition~\ref{prop:sym}, the derivative of the first two
terms of \eqref{eq:no gradient var} is
\begin{equation}\label{eq:derivative1}
  \begin{aligned}
    &\int_{\mathrlap{\Omega \backslash (A \cup B)}}\grad_\theta   q_\theta(\x)
    \left((I-P^i)q_\theta (\x)\right) \rho(\x) d \x -
    \int_{\mathrlap{\Omega \backslash (A \cup B)}} \grad_\theta   q_\theta(\x)
    P^b r(\x)\rho(\x)                 d \x\\
    =& \int_{\mathrlap{\Omega \backslash (A \cup B)}}\grad_\theta    q_\theta(\x)
    \left(
    q_\theta(\x)-\E^{\x}\left(q_\theta(\x_{\delta})\mathbf{1}_{\{\delta<\tau\}}\right)
    - \E^{\x}\left(g (\x_{\tau})\mathbf{1}_{\{\delta\geq\tau\}}\right)
    \right) \rho(\x) d
    \x.
  \end{aligned}
\end{equation}
Notice that $\eqref{eq:derivative1}$ is an integral with measure $\rho(\x) d\x$. Therefore, if
$\x\sim \rho$ and $\x_\delta$ follows \eqref{eq:trajectory}, \eqref{eq:derivative1} can be written
as
\begin{equation}\label{eq:MCint}
  \E_{\x\sim \rho}\grad_\theta
  q_\theta(\x)
  \left(
  q_\theta(\x)-\E^{\x}\left(q_\theta\left(\x_{\delta}\right)\mathbf{1}_{\{\delta<\tau\}}\right)
  - \E^{\x}\left(r(\x_{\tau})\mathbf{1}_{\{\delta\geq\tau\}}\right)
  \right).
\end{equation}
An unbiased estimator for \eqref{eq:derivative1} in an SGD method is therefore
\begin{equation}\label{eq:unbiased}
  \grad_\theta q_\theta(\x) \left(q_\theta(\x) - q_\theta(\x_\delta)\mathbf{1}_{\{\delta<\tau\}} - r(\x_{\tau})\mathbf{1}_{\{\delta\geq\tau\}}\right)
  = \grad_\theta q_\theta(\x) \left(q_\theta(\x) - q_\theta(\x_\delta)\mathbf{1}_{\{\delta<\tau\}} - \mathbf{1}_{\{\delta\geq\tau = \tau_B\}}\right).
\end{equation}
The significance of \eqref{eq:unbiased} is that one only needs to take the gradient once with
respect to the parameter $\theta$, whereas in \citep{khoo2019solving} one needs to take the gradient
of both $\x$ and $\theta$.

Let us comment on the implementation details of \eqref{eq:unbiased}. First, the sample $\x$ is
supposed to be sampled from $\rho$. As mentioned previously, if the potential function $V$ is
confining, then the Langevin process \eqref{langevin} is ergodic. This implies that the distribution
of the samples $\x_t$ generated from the stochastic differential equation (SDE)
\eqref{eq:trajectory} converges to $\rho$ in the limit of $t\rightarrow\infty$. Therefore, after an
initial burn-in period, the samples from the SDE trajectory follow the distribution $\rho$. There
has been a large literature on how to solve SDE \eqref{eq:trajectory} numerically. The simplest
method is the Euler-Maruyama scheme (see for example \citep{kloeden2013numerical}). Let $\Delta t>0$
be the time step size and $\mathbf{w}_{\Delta t}\sim\mathcal{N}(0, {\Delta t}I_d)$, and $I_d$ is the
$d$-dimensional identity matrix. The Euler-Maruyama approximation $\tilde{\x}_{n\Delta t}$ at time
steps $n\Delta t$ are computed via
\[
\tilde{\x}_{(n+1)\Delta t} = \tilde{\x}_{n\Delta t} -
\grad{V}(\tilde{\x}_{n\Delta t})\Delta t+\sqrt{2\beta^{-1}}\mathbf{w}_{\Delta t}.
\]
Assuming ergodicity, the distribution $\x\sim \rho$ can be approximated by
$\tilde{\x}_{N\Delta t}$ for a sufficiently small $\Delta t$ and sufficiently large $N$
with an arbitrary $\tilde{\x}_{0}$.

Second, given $\x\sim\rho$, one needs to sample $\x_\delta$. We approximate $\x_\delta$ by Euler-Maruyama scheme 
as well
\begin{equation}\label{eq:nextstep}
  \x_\delta = \x -\grad{V}(\x)\delta+\sqrt{2\beta^{-1}}\mathbf{w}_\delta, 
\end{equation}
where $\mathbf{w}_\delta\sim\mathcal{N}(0,\delta I_d)$.

Finally, it is necessary to determine the indicators $\mathbf{1}_{\{\delta<\tau\}}$,
$\mathbf{1}_{\{\delta\geq\tau = \tau_A\}}$ and $\mathbf{1}_{\{\delta\geq\tau = \tau_B\}}$ in order
to compute \eqref{eq:unbiased}. For this, the following approximations are used
\[
\mathbf{1}_{\{\delta<\tau\}} = 1 \;\text{if}\; \x_\delta\in\Omega\backslash A\cup B,\quad
\mathbf{1}_{\{\delta\geq\tau = \tau_A\}} = 1 \;\text{if}\; \x_\delta\in A, \quad
\mathbf{1}_{\{\delta\geq\tau = \tau_B\}} = 1 \;\text{if}\; \x_\delta\in B.
\]

Intuitively these approximations should work well when $\delta$ is sufficiently small, which is also
supported by the numerical experiments. Even for a fixed $\delta$, we can apply Euler-Maruyama
  scheme with multiple steps to improve the accuracy. Specifically, given $\mathbf{x}\sim\rho$, one
  can approximate $\mathbf{x}_\delta$ by the last term of the sequence $\x_0, \x_{\delta/M}, \ldots,
  \x_\delta$, where
\begin{equation}\label{eq:multiEM}
\mathbf{x}_{(k+1)\delta/M} = \x_{k\delta/M} -\grad{V}(\x_{k\delta/ M})\delta/M+\sqrt{2\beta^{-1}}\mathbf{w}_{\delta/M}, 
\end{equation}
and $\x_0 = \x$. When using the multiple-step approximation, the indicators are determined by 
\[
\begin{aligned}
\mathbf{1}_{\{\delta<\tau\}} &= 1 \;\text{if}\; \{\x_{k\delta/M}\}_{k=1}^{M}\subset\Omega\backslash A\cup B,\\
\mathbf{1}_{\{\delta\geq\tau = \tau_A\}} &= 1 \;\text{if}\; \{\x_{k\delta/M}\}_{k=1}^{M}\cap A\not=\emptyset,\\
\mathbf{1}_{\{\delta\geq\tau = \tau_B\}} &= 1 \;\text{if}\; \{\x_{k\delta/M}\}_{k=1}^{M}\cap B\not=\emptyset.
\end{aligned}
\]

It is worth mentioning that various importance sampling methods (see for example
\citep{li2019computing,rotskoff2020learning}) can potentially be used to approximate
the integral in \eqref{eq:MCint} and the corresponding gradients calculated above. One can generate the initial state $\x$ in \eqref{eq:nextstep} using importance sampling techniques, and then generate $\x_\delta$ and also subsequent samples according to \eqref{eq:nextstep}. 

\paragraph{Derivative of the penalty terms.}
For the third and fourth terms of \eqref{eq:no gradient var}, unbiased estimators of their gradients
are
\begin{equation}
  c \grad_\theta q_\theta(\x_A) q_\theta(\x_A),\quad
  c \grad_\theta q_\theta(\x_B) (q_\theta(\x_B)-1),
\end{equation}
respectively, where $\x_A\sim m_A$ and $\x_B\sim m_B$.

\paragraph{Connection with reinforcement learning.}
Solving committor function can be viewed as a special case of the policy evaluation problem in
reinforcement learning (RL): $\Omega \backslash (A \cup B)$ is the state space; the transition kernel
is given by the operator $P^i$ defined in \eqref{eq:pi}; the discount factor is equal to one;
there is no immediate reward for each individual step but the final reward is $(P^b g)(\x)$;
$q(\x)$ is the value function of the problem; \eqref{eq:nobd} is the Bellman equation; and
\eqref{eq:unbiased} is the temporal difference (TD) update.

However, the key difference with the general policy evaluation problem is that, due to the detailed
balancing of overdamped Langevin dynamics, there is a variational formulation \eqref{eq:variational}
for \eqref{eq:nobd}, which is not available for general RL problems except the special case
considered in \citep{ollivier2018approximate}. Because of the variational formulation,
\eqref{eq:unbiased} is guaranteed to converge at least to a local minimum, even in the neural
network parameterization. On the other hand, such a stability result fails to exist for neural
network approximations for general RL problems, as many techniques and tricks are required in order
for the neural network parameters to converge.

\subsection{Neural network architecture}\label{sec:arch}
\label{sec:achitecture}

The architecture of the neural network $q_\theta(\x)$ follows the one used in
\citep{khoo2019solving} and is specifically designed for this problem. Below we briefly summarize
the main design considerations.

\begin{figure}[ht]
  \centering
  \includegraphics[width=0.8\textwidth]{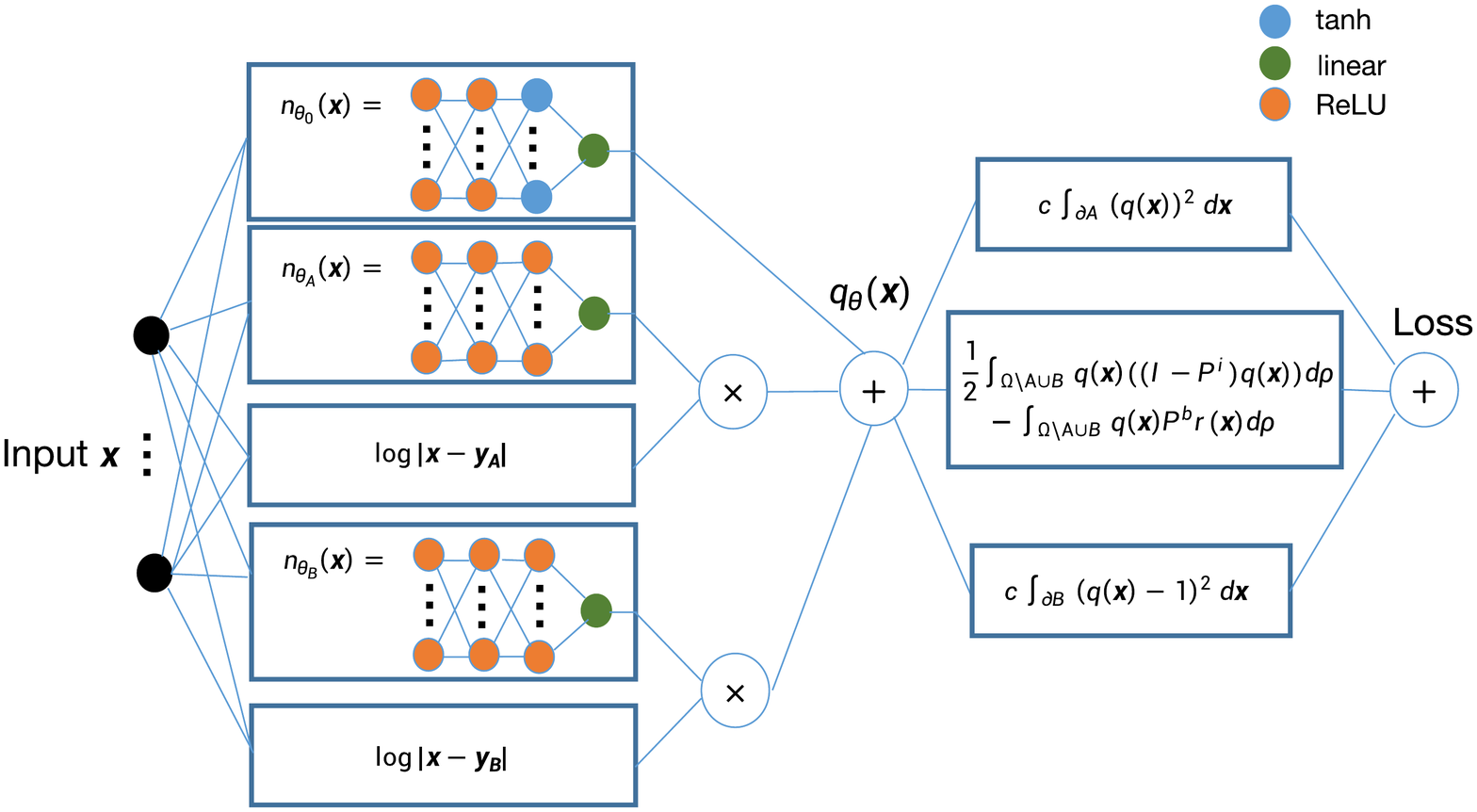}
  \caption{An example of the neural network architecture for a committor function and the
    corresponding loss. In this example we have two $\log|\x-\mathbf{y}|$ type
    singularities.}
  \label{fig:nn}
\end{figure}

In the low temperature regime, i.e. when $T\rightarrow 0$, there is typically a sharp
interface between $A$ and $B$, as pointed out by \citep{khoo2019solving}. In order to account for
this sharp transition, the tangent function $\tanh$ is used as the activation function at the last
nonlinear layer, as shown in Fig.~\ref{fig:nn}.

In the high temperature regime, i.e. when $T\rightarrow \infty$, there is another type of singular
behavior. As $T\rightarrow\infty$, the equation \eqref{committor} converges heuristically to a
Laplace equation with a Dirichlet boundary condition. When the domains $A$ and $B$ are relatively
small, the solution near the boundary $\partial A$ and $\partial B$ are dictated asymptotically by
the fundamental solution
\begin{equation}\label{eq:fund}
  \Phi(\x):=\left\{\begin{array}{ll}
  -\frac{1}{2 \pi} \log |\x| & (n=2) \\
  \frac{\Gamma(n/ 2)}{(2 \pi)^{n / 2}|\x|^{n-2}} & (n \geq 3)
  \end{array}\right.,
\end{equation}
where $n$ is the inherent dimension of the Laplace equation of $T\rightarrow\infty$, as demonstrated
in \citep{khoo2019solving}. For example, in the rugged-Muller problem \ref{sec:RM}, the inherent
dimension $n=2$, while in the Ginzburg-Landau problem \ref{sec:G-L} we have $n=d=49$.

In order to address these two types of singular behaviors, we introduce the following
parameterization
\begin{equation}\label{eq:NNpara}
  q_{\theta}(\x):=
  n_{\theta_{A}}(\x) S_A\left(\x-\mathbf{y}^{A}\right) +
  n_{\theta_{B}}(\x) S_B\left(\x-\mathbf{y}^{B}\right) + n_{\theta_{0}}(\x),
\end{equation}
where $\mathbf{y}^A$ and $\mathbf{y}^B$ are the centers of $A$ and $B$,  $S_A(\x-\mathbf{y})$ and
$S_B(\x-\mathbf{y})$ are set to be fundamental solutions \eqref{eq:fund}, with $n$ depending on the
inherent dimension of the problem. Finally, $n_{\theta_{A}}$, $n_{\theta_{B}}$ are fully connected
neural networks with ReLU activation, and $n_{\theta_{0}}$ is a fully connected neural network with
$\tanh$ activation at the last nonlinear layer and ReLU activation at other nonlinear layers. This
architecture is summarized in Fig.~\ref{fig:nn}.

\section{Lazy training analysis of the optimization}\label{sec:lazy}

When the learning rate approaches zero, the dynamics of the SGD can be approximated by the
corresponding gradient flow (see \citep{kushner2003stochastic} for example):
\begin{equation}\label{eq:vanillagf}
    \theta'(t) = -\grad_\theta R(q_\theta), \quad \theta(0) = \theta_0, 
\end{equation}
where in our case the loss function $R(\cdot)$ takes the form
\begin{equation}\label{eq:R}
  R(q) = \frac{1}{2} \langle (I-P^i)q-P^b r, q\rangle_{\rho}
  +\frac{c}{2} \langle q, q\rangle_{m_A}
  +\frac{c}{2} \langle q-1, q-1\rangle_{m_B},
\end{equation}
with $\langle u, v\rangle_{m_A} = \int uvdm_A$ and $\langle u, v\rangle_{m_B} = \int uvdm_B$. In
this section, we consider instead the rescaled gradient flow
\begin{equation}\label{eq:scaledgf}
  \theta'(t) = -\frac{1}{\alpha^2}\grad_{\theta} R(\alpha q_{\theta}), \quad \theta(0) = \theta_0, 
\end{equation}
with $ R(\alpha q) = \frac{1}{2} \langle (I-P^i)(\alpha q)-P^b r, \alpha q\rangle_{\rho} +\frac{c}{2}
\langle \alpha q, \alpha q\rangle_{m_A} +\frac{c}{2} \langle \alpha q-1, \alpha q-1\rangle_{m_B} $, 
and analyze the training dynamics of $q_\theta$ in the lazy training regime. The reason for
considering this rescaled formula is that the scaling effect caused by $\alpha$ arises in several
situations, for example, when the weights of the NN are large in magnitude at initialization and the
learning rate is small (see for example, \citep{chizat2019lazy,agazzi2020temporaldifference}).

In \citep{chizat2019lazy}, it has been shown that when the scaling factor $\alpha$ is sufficiently
large, the gradient flow \eqref{eq:scaledgf} converges at a geometric rate to a local minimum of
$F_\alpha(\theta) \equiv R(\alpha q_{\theta})/\alpha^2$, under some conditions that are detailed
below.

\begin{theorem}\citep{chizat2019lazy}\label{thm:lazy}
  Assume that: (1) $\theta\mapsto q_\theta \in \mathcal{F}$, where $\mathcal{F}$ is a separable Hilbert space; (2) $q$ is differentiable with a locally Lipschitz differential $\grad_\theta q_\theta$; (3) $q_{\theta_0} =
  0$ and $\operatorname{rank} \grad_\theta q_\theta$ is a constant in a neighborhood of $\theta_0$; (4) $R$ is strongly
  convex and differentiable with a Lipschitz gradient.

Then there exists $\alpha_0 > 0$, such that for any $\alpha>\alpha_0$, the gradient flow \eqref{eq:scaledgf} converges at a geometric rate to a local minimum of $F_\alpha$.
\end{theorem}

In our setting, $\mathcal{F}$ is the separable Hilbert space $L_{\nu}^2((\Omega\backslash (A\cup
B))\cup \partial A \cup \partial B)$, where $\nu = \rho+cm_A+cm_B$ and the rescaled gradient flow is
given by
\begin{equation}\label{eq:gradflow}
\begin{aligned}
  \theta'(t) =-\frac{1}{\alpha}\Big[&\langle (I-P^i)\alpha q_{\theta(t)}-P^b r, \grad_\theta q_{\theta(t)}\rangle_{\rho}
    -c \langle \alpha q_{\theta(t)}, \grad_{\theta} q_{\theta(t)}\rangle_{m_A}
    -c \langle \alpha q_{\theta(t)}-1, \grad_\theta q_{\theta(t)}\rangle_{m_B}\Big]. 
    \end{aligned}
\end{equation}
The following result states the strong convexity of $R$.
\begin{prop}\label{prop:strcvx}
  Assume that the operator $P^i$ defined in \eqref{eq:pi} has a corresponding probability density
  function $p^i(\x, \mathbf{y})$ such that
  \begin{equation}\label{eq:pdf}
    P^{i}f(\x) = \int_{\mathrlap{\Omega\backslash (A\cup B)}} p^i(\x, \mathbf{y})f(\mathbf{y})d\mathbf{y}, 
  \end{equation}
  and 
  \begin{equation}\label{eq:HS}
    \iint_{\mathrlap{\Omega\backslash (A\cup B)\times\Omega\backslash (A\cup B)}} p^i(\x, \mathbf{y})p^i(\mathbf{y}, \x)d\x d\mathbf{y}<\infty.
  \end{equation}
  Then the quadratic form $f\mapsto \langle f, (I-P^i)f\rangle_{\rho}$ is strongly convex on $L_{\rho}^2(\Omega\backslash A\cup B)$. 
\end{prop}

Based on Proposition~\ref{prop:strcvx}, the following theorem is a direct consequence of
Theorem~\ref{thm:lazy}.
\begin{theorem}
Assume that \eqref{eq:pdf} and \eqref{eq:HS} hold, and $q_{\theta_0} = 0$, and that
$\operatorname{rank} \grad_\theta q_\theta$ is a constant in a neighborhood of $\theta_0$, then there exists
$\alpha_0>0$, such that for any $\alpha>\alpha_0$, the gradient flow \eqref{eq:gradflow}

converges at a geometric rate to a local minimum of $F_\alpha(\theta) = R(\alpha q_\theta)/\alpha^{2}$.
\end{theorem}

\section{Numerical experiments} \label{sec:num}

This section presents the numerical results of the proposed method on several examples. The neural
network architecture follows Fig.~\ref{fig:nn} and the stochastic gradient is computed following the
description in Section~\ref{sec:nonlinear}. The training is carried out with the Adam optimizer
\citep{kingma2014adam}. Whenever the true solution $q^*$ is available, the performance is evaluated
by the relative error metric:
\begin{equation}\label{eq:error}
  E=\frac{\|q_{\theta}-q^*\|_{L_\rho^2(\Omega\backslash A\cup B)}}{\|q^*\|_{L_\rho^2(\Omega\backslash A\cup B)}}
\end{equation}
computed on a validation dataset generated according to the distribution $\rho$. In all the
experiments, $2000$ samples are used for the boundary measures $m_A$ and
$m_B$.

\subsection{The double-well potential}\label{sec:wwell}
In this experiment, consider the committor function in the following double-well potential:
\begin{equation}
  V(\x)=\left(x_{1}^{2}-1\right)^{2}+0.3 \sum_{i=2}^{d} x_{i}^{2}, 
\end{equation}
with $d=10$ and the regions $A$ and $B$ defined as
\begin{equation}
  A=\left\{x \in \mathbb{R}^{d} \mid x_{1} \leq-1\right\}, \quad B=\left\{x \in \mathbb{R}^{d} \mid
  x_{1} \geq 1\right\}.
\end{equation}
The true solution can be easily obtained by letting $q(\x) = \tilde{q}(x_1)$, such that it
satisfies the one-dimensional ordinary differential equation (ODE)
\begin{equation}\label{eq:ode}
  \frac{d^2 \tilde{q}(x_1)}{dx_1^2} - 4\beta x_1(x_1^2-1)\frac{d \tilde{q}(x_1)}{d x_1} = 0, \quad \tilde{q}(0)=0, \quad \tilde{q}(1) = 1.
\end{equation}
By solving this ODE numerically, we obtain a highly accurate approximation of the exact solution
$q^*$.

Since the solution of this problem does not exhibit any singular behavior at the boundaries $\partial A$ and $\partial B$, $S_A(\cdot) = S_B(\cdot)
= 0$ and the training is only performed for the component $n_{\theta_{0}}(\cdot)$ in
\eqref{eq:NNpara}. The problem is solved for the temperatures $T=0.5$ and
$T=0.2$. Table~\ref{table:Wwell} summarizes the numerical results and the training
parameters. Notice that when the temperature $T$ is lower, the distribution $\rho$ is sparser in $\Omega \backslash (A \cup B)$. Therefore, more samples are used for $T=0.2$ than for $T=0.5$ in
order to achieve a comparable precision. For example, for $T=0.5$, $1.5\times 10^5$ samples are used
for $\Omega \backslash (A \cup B)$ and the relative error is $E=0.014$. On the other hand, for
$T=0.2$, $8.0\times 10^5$ samples are used and the relative error of the neural network solution is
$E =0.011$.

\begin{table}[ht]
  \centering 
  \begin{tabular}{ c c c c c c c } 
	\hline\hline
	$T$  & $E$ & $\delta$ & $c$ & No. training samples & Batch size & No. testing samples  \\
	\hline\hline 
	0.5 & 0.014 & $0.003$ & 15 & $1.5\times 10^5$ & $1000$ & $4.0\times 10^5$\\
	\hline
	0.2 & 0.011 & $0.003$ & 15 & $8.0\times 10^5$ & $1000$ & $8.0\times 10^5$\\
	\hline
  \end{tabular}
  \caption{Results for the double-well potential problem.}\label{table:Wwell}
\end{table}

\begin{figure}[ht]
  \centering
  \includegraphics[width=0.4\textwidth]{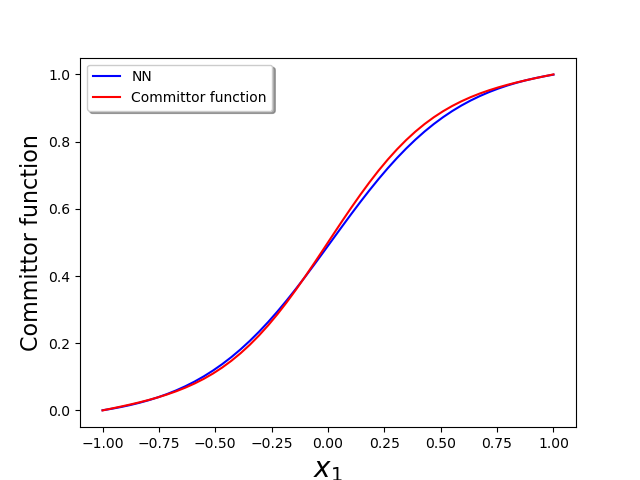}
  \caption{The committor function for the double-well potential along $x_1$ dimension when $T = 0.5$
    for an arbitrarily chosen $(x_2,\ldots,x_d)$ with $d = 10$.}
  \label{fig:Wwell}
\end{figure}

The numerical solution for $T=0.5$ is shown in Fig.~\ref{fig:Wwell}. The committor function
$q_\theta$ represented by a neural network and the committor function $q^*$ obtained via solving
\eqref{eq:ode} are plotted along the $x_1$ dimension for a fixed $(x_2,\ldots,x_d)$. The plot
demonstrates that the NN committor function gives a satisfactory approximation to the true
solution. We comment here that the final error is not sensitive to the parameter $\delta$. For
  example, for the case with $T=0.5$, $\delta = 0.003$ is used. If $\delta=0.01,0.03,0.05$ are
  chosen instead, the corresponding final errors are $E=0.013,0.013,0.013$, respectively. We also
  comment that using the multi-step sampling method \eqref{eq:multiEM} can slightly improve the
  final accuracy. For example, for the case with $T=0.5$ and $\delta = 0.05$, if the multiple-step
  scheme \eqref{eq:multiEM} with $M=10$ is used, the final error is $E=0.012$.

\subsubsection{Comparison with the method of \citep{khoo2019solving}}\label{sec:compare}
In this section, we compare the proposed method with the one in \citep{khoo2019solving} in
  terms of speed, accuracy, and robustness on the double-well potential problem with $T=0.5$. We use
  the same number of training samples ($1.5\times 10^5$) and testing samples ($4.0\times 10^5$), NN
  architectures, and hyperparameters in the implementation of both methods. The numerical tests are
  carried out on $4$ N$1$ virtual CPUs on the Google Cloud platform with altogether $26$ GB memory
  and a Tesla K80 GPU. In order to compare penalty coefficients on the same scale, a normalized
  penalty coefficient $c_{norm}$ is used: it is defined via $c_{norm}=c/\delta$ for the new method
  with $c$ given in \eqref{eq:no gradient var} and via $c_{norm}=\tilde{c}$ for the old method with
  $\tilde{c}$ given in \eqref{weakpenalty}.

\begin{figure}[ht]
  \centering \subfigure[Training loss of the proposed method\label{fig:compare_new}]{
    \includegraphics[width=0.46\textwidth]{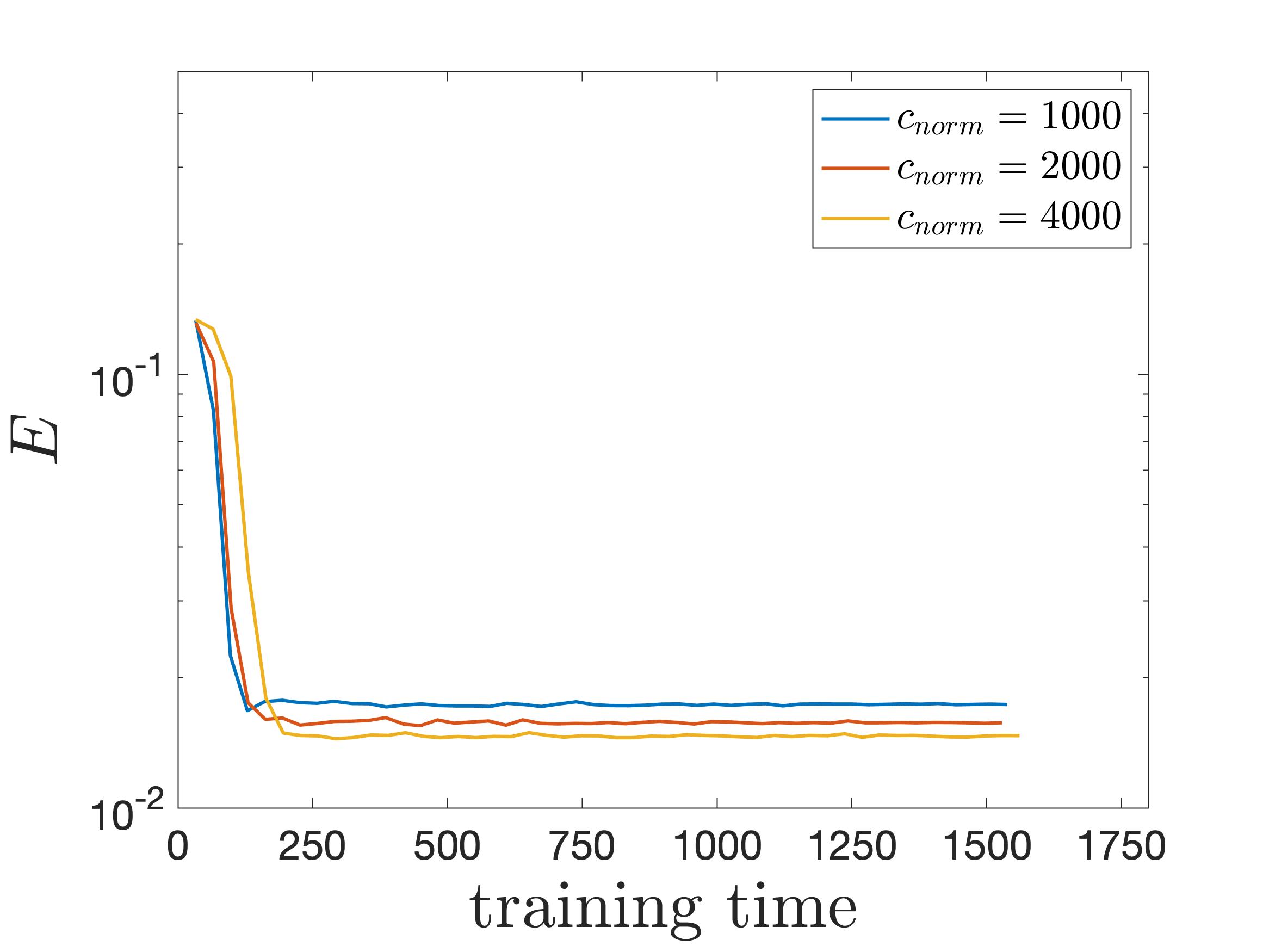} }
  \hspace{0.7em} \subfigure[Training loss of the method in \citep{khoo2019solving}
  \label{fig:compare_old}]{
    \includegraphics[width=0.46\textwidth]{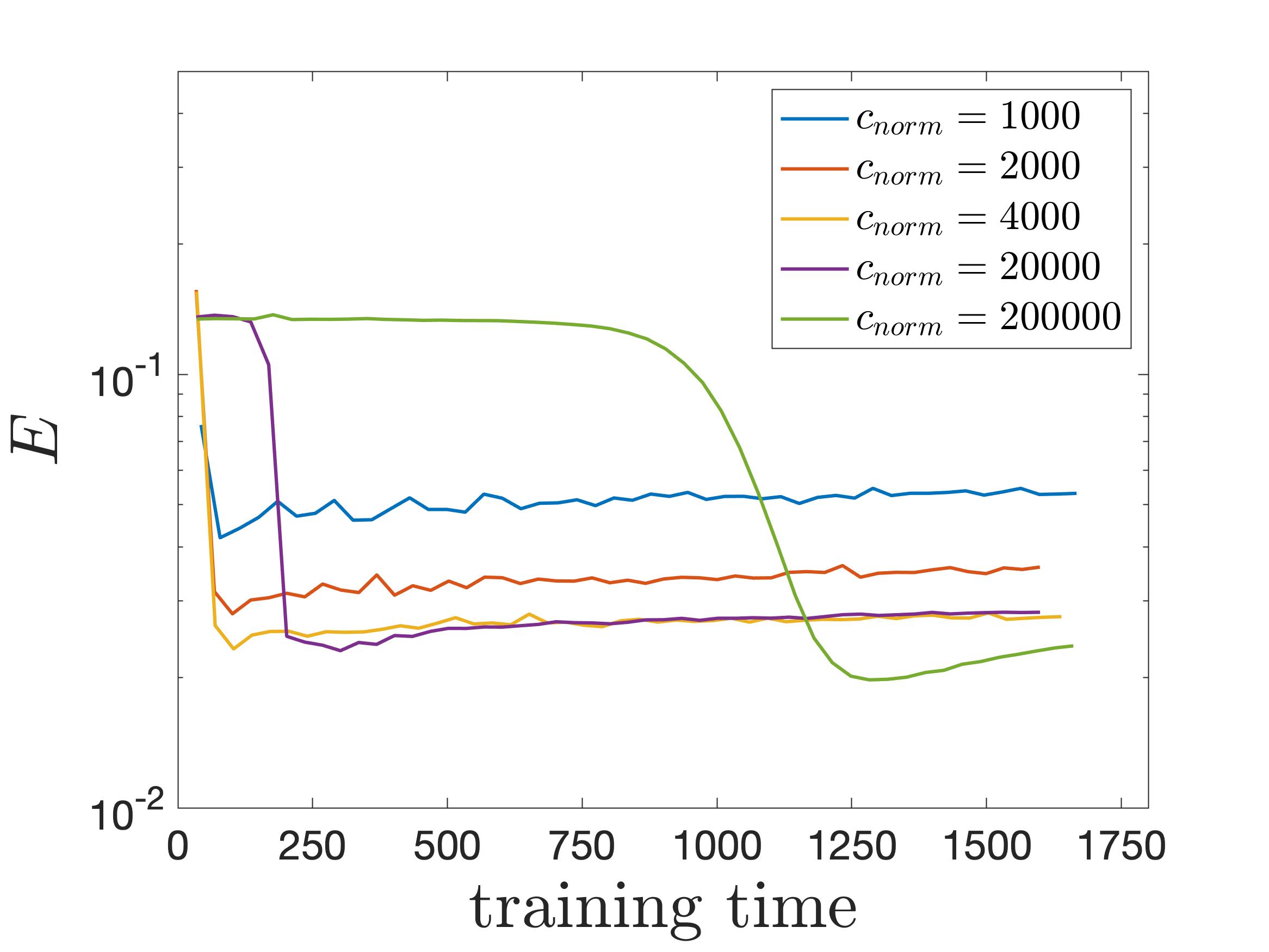} }
  \caption{Comparison of the training process of the proposed method and \citep{khoo2019solving}'s method. 
  (a): The training loss versus the training time used for the proposed method. 
  (b): The training loss versus the training time used for the method in \citep{khoo2019solving}. 
  We report the error $E$ (defined in \eqref{eq:error}) versus the training time used. 
  Here, $c_{norm}$ stands for the normalized penalty coefficient defined in Section~\ref{sec:compare}. }
  \label{fig:compare}
\end{figure}

 Fig.~\ref{fig:compare} demonstrates a clear difference between the behavior of the two methods
  during the training. As shown in Fig.~\ref{fig:compare_new}, when using the method proposed in
  this paper, the approximate solution converges quickly and the final relative error is rather
  small, regardless of the choice of penalty coefficients. In contrast, as shown in
  Fig.~\ref{fig:compare_old}, when using the method proposed in \citep{khoo2019solving}, different
  penalty coefficients lead to different training behaviors. When the penalty parameter $c_{norm}$
  is small, the time used to reach convergence is short but the final relative error is relatively
  large. When a large penalty parameter is used, the relative error is reduced but it is still
  higher than the proposed method. Moreover, the time for training to converge is long when using a
  large penalty parameter. In conclusion, when using the method in \citep{khoo2019solving}, the
  penalty coefficient needs to be carefully tuned in order to have a performance close to the
  proposed method.

\subsection{The rugged-Muller potential}\label{sec:RM}

In this example, we consider the committor function corresponding to the following rugged-Muller
potential:

\begin{equation}\label{eq:rmpotential}
  \begin{aligned}
    V(\x)=\tilde{V}(x_1, x_2)+\frac{1}{2 \sigma^{2}} \sum_{i=3}^{d} x_{i}^{2}, 
  \end{aligned}
\end{equation}
where 
\begin{equation}
  \tilde{V}(x_1, x_2) =
  \sum_{i=1}^{4} D_{i} e^{
    a_{i}\left(x_{1}-X_{i}\right)^{2}+b_{i}\left(x_{1}-X_{i}\right)\left(x_{2}-Y_{i}\right)
  +c_{i}\left(x_{2}-Y_{i}\right)^{2}} + \gamma\sin (2 k \pi x_1) \sin (2 k \pi x_2)
\end{equation}
is the $2$-dimensional rugged-Muller potential with the parameters
\begin{equation}
  \begin{aligned}
    \left[a_{1}, a_{2}, a_{3}, a_{4}\right] &=[-1,-1,-6.5,0.7], & 
    \left[b_{1}, b_{2}, b_{3}, b_{4}\right] &=[0,0,11,0.6], \\
    \left[c_{1}, c_{2}, c_{3}, c_{4}\right] &=[-10,-10,-6.5,0.7], &
    \left[D_{1}, D_{2}, D_{3}, D_{4}\right] &=[-200,-100,-170,15], \\
    \left[X_{1}, X_{2}, X_{3}, X_{4}\right] &=[1,0,-0.5,-1], &
    \left[Y_{1}, Y_{2}, Y_{3}, Y_{4}\right] &=[0,0.5,1.5,1],\\
    \left[\gamma,k,\sigma, d\right] &=[9,5,0.05, 10].
  \end{aligned}
\end{equation}
The domain of interest $\Omega$ of this example is $[-1.5, 1]\times[-0.5, 2]\times\mathbb{R}^{d-2}$
and the regions $A$ and $B$ are the following two cylinders:
\begin{equation}\label{eq:AB}
  \begin{aligned}
    A&=\left\{\x \in \mathbb{R}^{d} \mid \sqrt{(x_{1}+0.57)^2+(x_{2}-1.43)^2}\leq 0.3\right\}, \\
    B&=\left\{\x \in \mathbb{R}^{d} \mid \sqrt{(x_{1}-0.56)^2+(x_{2}-0.044)^2}\leq 0.3\right\}.
\end{aligned}
\end{equation}

In order to compute the error, $q^*$ is solved approximately within the $x_1x_2$-plane. More
precisely, we first apply finite element method on uniform grid to \eqref{committor} in $2$
dimensions with the potential $\tilde{V}$, the domain $\tilde{\Omega} = [-1.5, 1]\times[-0.5, 2]$,
and the regions $\tilde{A}$ and $\tilde{B}$ being the projection of the $A$ and $B$ defined in
\eqref{eq:AB} onto the $x_1x_2$-plane. Once the $2$-dimensional committor function $\tilde{q}$ is
available, the approximation is $q^*(\x) = \tilde{q}(x_1, x_2)$. The code of the finite element
method is provided by the authors of \citep{lai2018point}.

As mentioned in Section~\ref{sec:arch}, the singularity functions $S_A$ and $S_B$ should be set as
the fundamental solutions of \eqref{committor} when taking $T\rightarrow\infty$. In this case, the
limiting committor function as $T\rightarrow\infty$ has the form $q(\x) = \tilde{q}(x_1,x_2)$ that
satisfies a $2$-dimensional Laplace equation, the fundamental solution to which has the form
\[
\Phi(x_1, x_2) = -\frac{1}{4\pi}\log((x_1-a)^2+(x_2-b)^2).
\]
Therefore, we set $S_A = \log ((x_{1}+0.57)^2+(x_{2}-1.43)^2)$ and
$S_B=\log((x_{1}-0.56)^2+(x_{2}-0.044)^2)$.

\begin{figure}[ht]
  \centering
  \subfigure[$T = 22$ committor function]{ \includegraphics[width=0.48\textwidth]{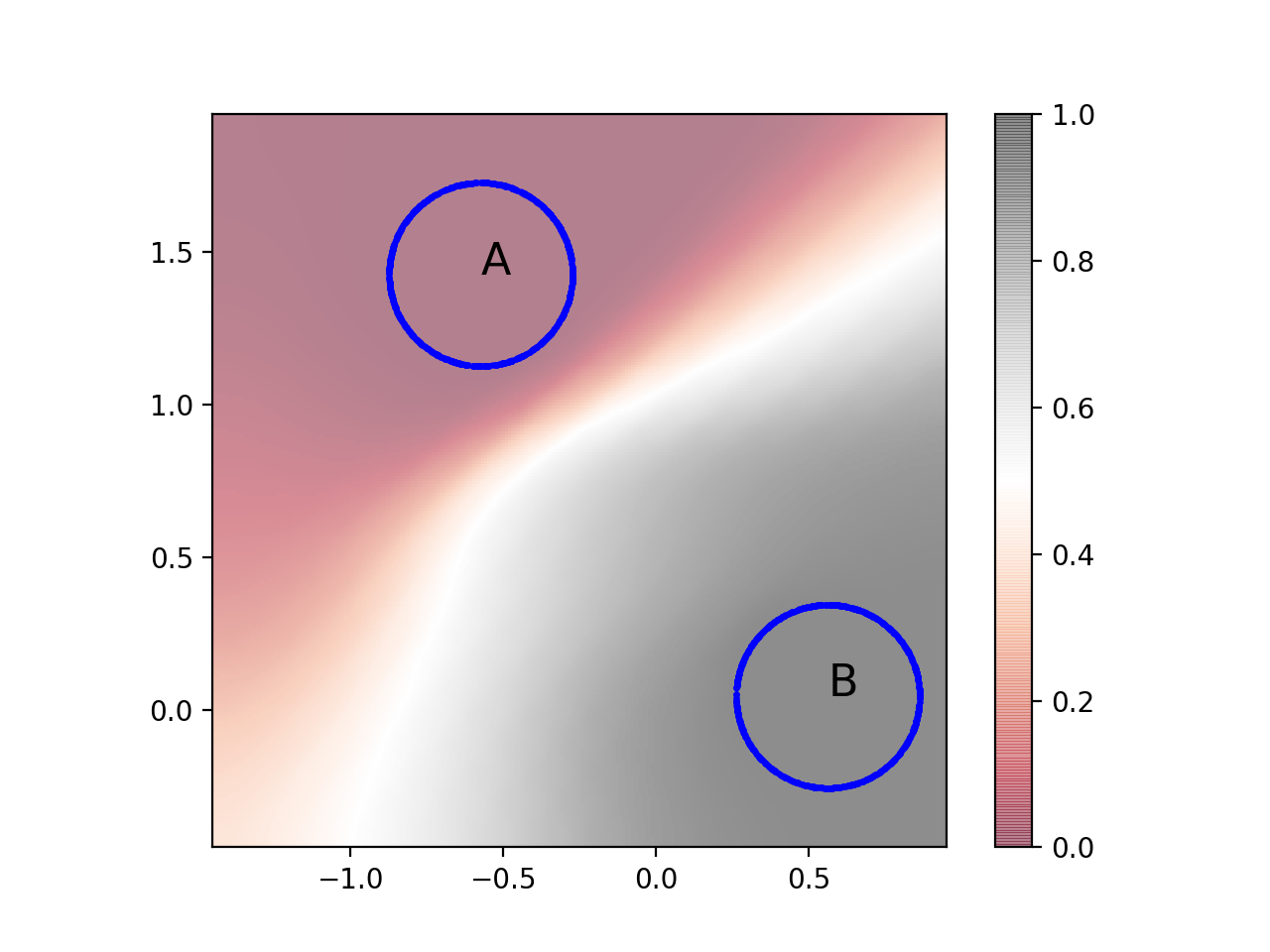} }\label{fig:rmtrue}
  \subfigure[$T = 22$ NN approximation]{ \includegraphics[width=0.48\textwidth]{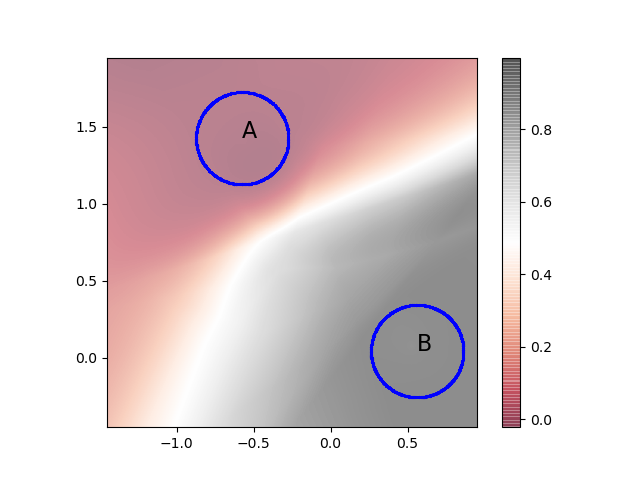} }\label{fig:rmnn}
  \subfigure[$T = 40$ committor function]{ \includegraphics[width=0.48\textwidth]{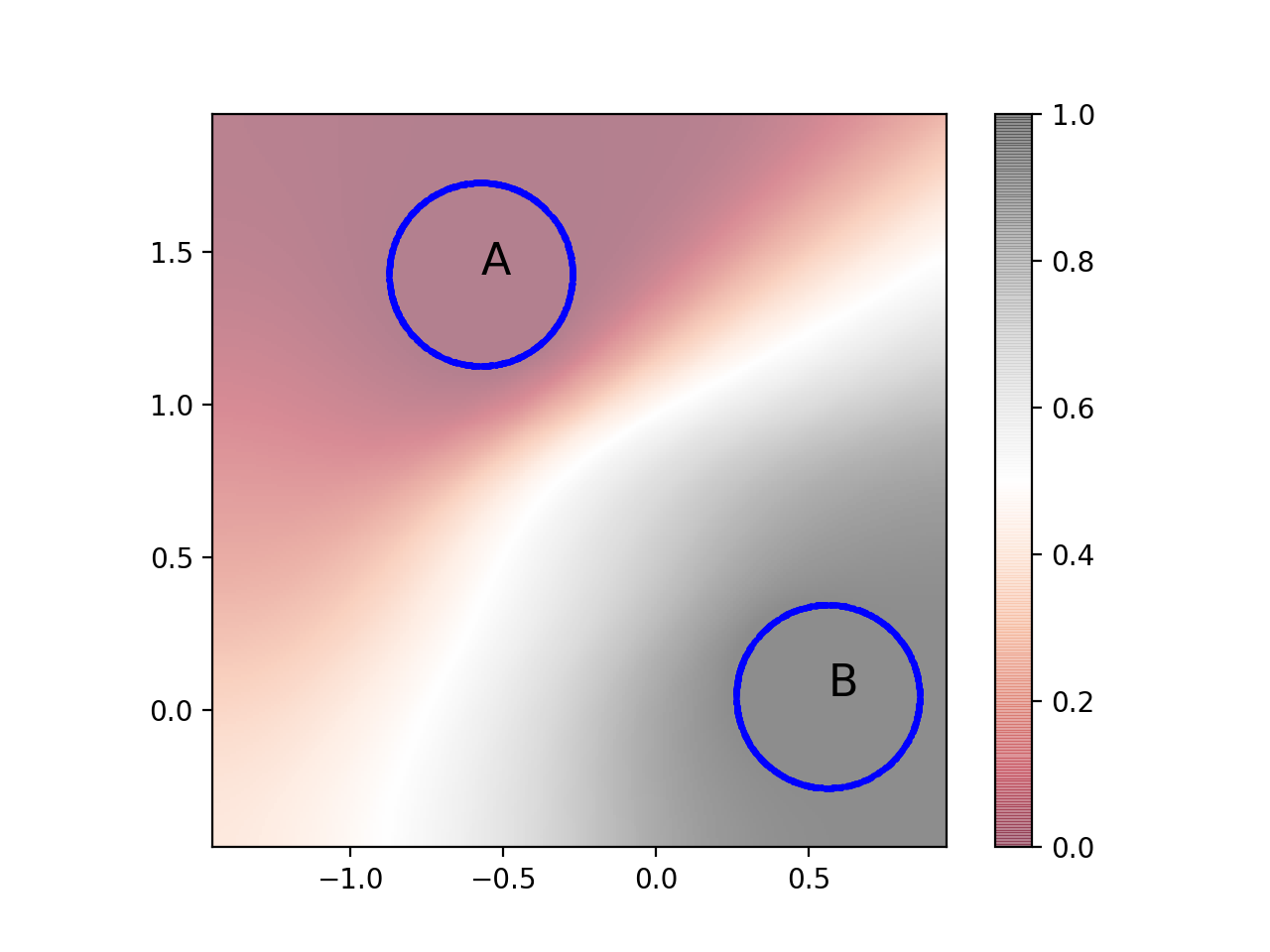} }\label{fig:rmtrue40}
  \subfigure[$T = 40$ NN approximation]{ \includegraphics[width=0.48\textwidth]{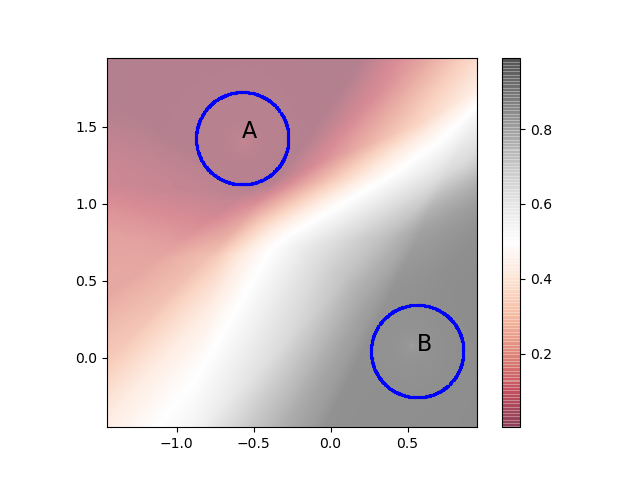} }\label{fig:rmnn40}
  \caption{\label{fig:rm}
    Comparisons between the NN represented committor functions and the ground truths.
    (a): The ground truth committor function for $T = 22$. (b): The NN parameterized committor
    function for $T = 22$. (c): The ground truth committor function for $T = 40$. (d): The NN
    parameterized committor function for $T = 40$.}
\end{figure}

Numerical experiments are carried out when $T=22$ and $T=40$. In both situations, $6.0\times 10^5$
samples are used in $\Omega \backslash (A \cup B)$. The relative errors are $E = 0.024$ when $T=22$,
and $E = 0.023$ when $T=40$. Table~\ref{table:RM} summarizes the numerical error and the parameters
used in the experiments.
\begin{table}[ht]
  \centering 
  \begin{tabular}{ c c c c c c c } 
	\hline\hline 
	$(T, \sigma)$  & $E$ & $\delta$ & $c$ & No. training samples & Batch size & No. testing samples \\
    \hline\hline 
	(22, 0.05) & 0.024 & $0.001$ & 500 & $6.0\times 10^5$ & $5000$ & $1.0\times 10^6$\\
	(40, 0.05) & 0.023 & $0.001$ & 500 & $6.0\times 10^5$ & $5000$ & $1.0\times 10^6$\\
	\hline 
  \end{tabular}
  \caption{Results for the rugged-Muller potential problem. }\label{table:RM}
\end{table}
The numerical solution is plotted in Fig.~\ref{fig:rm}, where the exact solutions $q^*$ are shown on
the left and the committor functions $q_\theta$ represented by neural network are plotted on the
right. The plots show that the NN approximation shows good agreement with the true solution.

\subsection{The Ginzburg-Landau model}\label{sec:G-L}

The Ginzburg-Landau theory is developed to give a mathematical description of superconductivity
\citep{hoffmann2012ginzburg}.  In this example, we discuss a simplified Ginzburg-Landau phase
transition model. The Ginzburg-Landau energy in one dimension is defined as:

\begin{equation}
  \label{eq:1D G-L}
  \tilde{V}[u]=\int_0^1 \dfrac{\lambda}{2} u_x^2+\dfrac{1}{4\lambda}(1-u^2)^2 d x, 
\end{equation}
where $\lambda$ is a small positive parameter and $u$ is a sufficiently smooth function on $[0,1]$
with boundary conditions $u(0) = u(1) = 0$.

\begin{figure}[!htbp]
  \centering
  \subfigure[Local minimizer $u_-$]{ \includegraphics[width=0.45\textwidth]{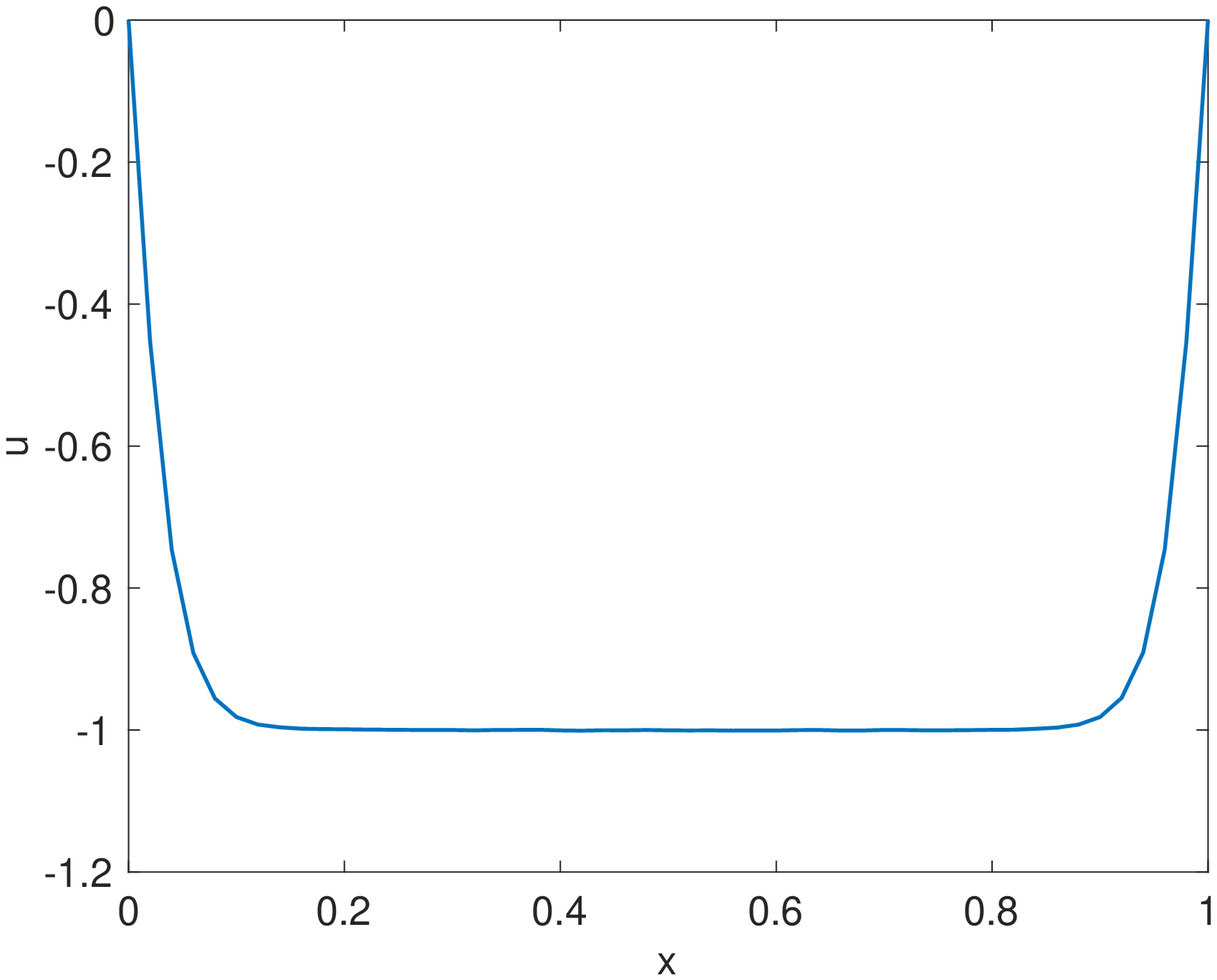} }\label{fig:minima-}
  \subfigure[Local minimizer $u_+$]{ \includegraphics[width=0.45\textwidth]{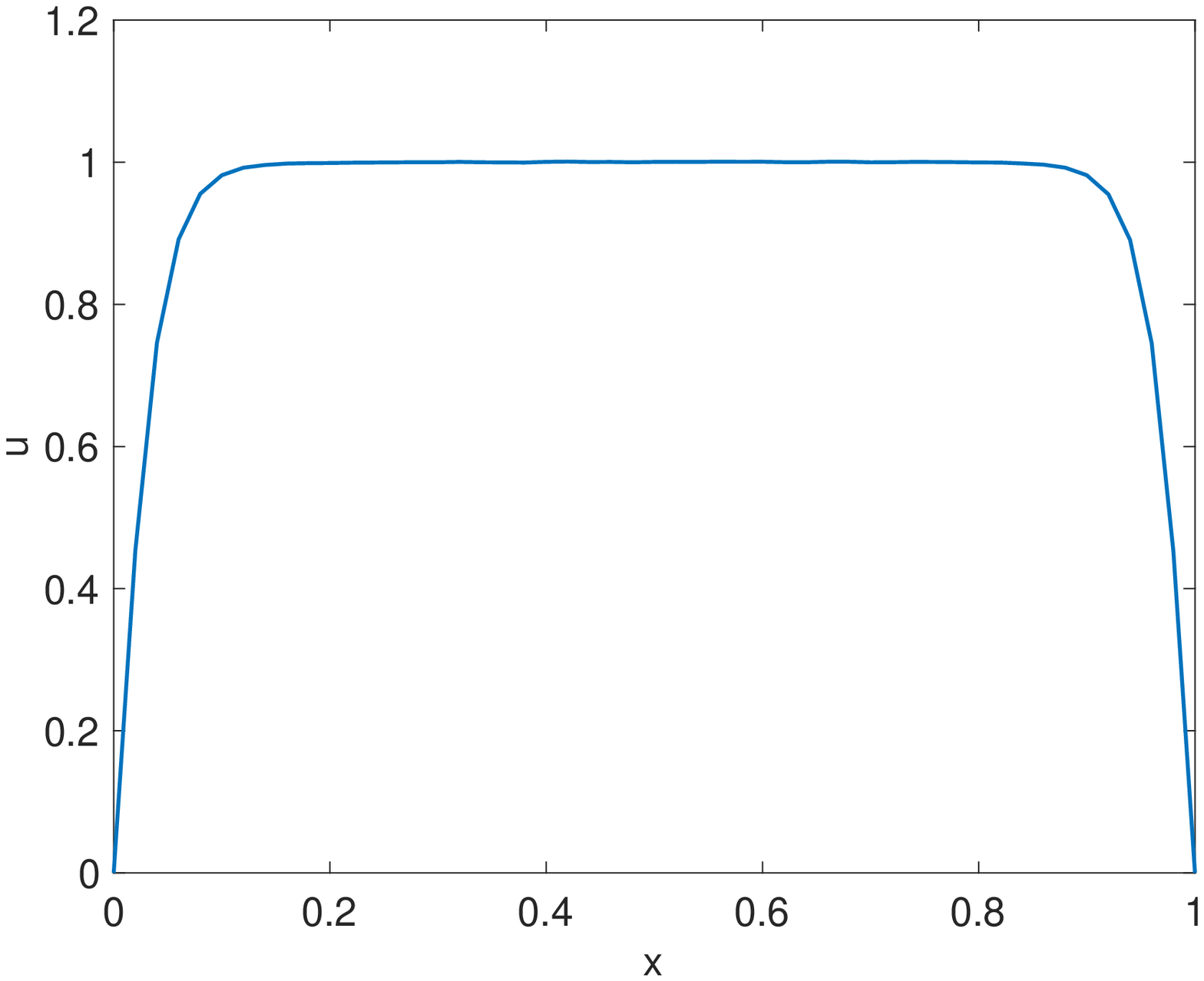} }\label{fig:minima+}
  \caption{\label{fig:minima} Two local minima of the energy \eqref{eq:1D discret G-L} with
    $\lambda=0.03$. (a): $u_-$, (b): $u_+$.  }
\end{figure}

The high-dimensionality nature of the committor functions is a direct result of the discretization
of $u$. With a numerical discretization, $u(x)$ is uniformly discretized by $U=(U_1,\cdots,U_d)$
defined on a uniform grid on $[0,1]$ with the boundary conditions $U_0=U_{d+1}=0$. Then the
continuous Ginzburg-Landau energy is approximated by a discrete one:
\begin{equation}
  \label{eq:1D discret G-L}
  V(U):=\tilde{V}_h[U]=\sum_{i=1}^{d+1} \dfrac{\lambda}{2}\left(\dfrac{U_i-U_{i-1}}{h}\right)^2
  +\dfrac{1}{4\lambda}(1-U_i^2)^2,
\end{equation}
where the grid size $h=1/(d+1)$. In this experiment we use $h = 1/50$ and the dimension
$d=49$. $V(U)$ has two local minima $u_\pm(\cdot)$ shown in Fig.~\ref{fig:minima}. The regions A and
B are taken as the spheres $\{U:||U-u_\pm||\leq r\}$, where $\| \cdot \|$ is the Euclidean norm, and
the radius $r$ is chosen to be $3$.

Based on the discretization used, $d=49$. According to Section~\ref{sec:achitecture}, the
singularities $S_A$ and $S_B$ are set to $S_A = |U-u_-|^{2-d}$, $S_B = |U-u_+|^{2-d}$.  The
numerical results are reported for the temperatures $T=30$ and $T=20$. For both cases $2.0 \times
10^5$ samples are used in $\Omega \backslash (A \cup B)$. In Table~\ref{table:GL} we summarize the
parameters used in the experiments with $T=30$ and $T=20$.

\begin{table}[ht]
  \centering 
  \begin{tabular}{ c c c c c c} 
	\hline\hline 
	$T$   & $\delta$ & $c$ & No. training samples & Batch size  \\ %
	\hline\hline 
	20  & $0.002$ & 200 & $2.0\times 10^5$ & $5000$ \\
	30  & $0.001$ & 200 & $2.0\times 10^5$ & $5000$ \\
	\hline 
  \end{tabular}
  \caption{Parameters for the Ginzburg-Landau problem.}\label{table:GL}
\end{table}

In this problem, it is intractable to obtain the exact $q^*$ due to the high dimensionality and
therefore we are not able to estimate the relative error $E$ directly. Instead, we study the region
near the $\frac{1}{2}$-isosurface of committor function $q_\theta$, which is defined as
$\Gamma_{\frac{1}{2},\epsilon} = \{U:\vert q_\theta(U)-\frac{1}{2}\vert < \epsilon\}$. If $q_\theta$ is indeed a
satisfactory approximation of $q^*$, then for a trajectory given by \eqref{eq:trajectory} starting
from an arbitrary point $\x_0\in\Gamma_{\frac{1}{2},\epsilon}$, the probability of entering $B$
before $A$ should be close to $\frac{1}{2}$.

More precisely, we first identify $m$ states $\{\tilde{\x}_j\}_{j=1}^m$ on
$\Gamma_{\frac{1}{2},\epsilon}$. From each $\tilde{\x}_j$, $N$ trajectories are generated
according to \eqref{eq:trajectory}. Let us denote the number of trajectories reaching $B$ before $A$
as $n$. If the NN committor function is accurate, then by the central limit theorem, when $N$ is
large, the distribution of $n/N$ should be approximately $\mathcal{N}(\frac{1}{2},(4N)^{-1})$,
i.e. the normal distribution with mean $\frac{1}{2}$ and variance $(4N)^{-1}$.  In the actual
experiment with $\epsilon=0.01$, $m=120$, and $N=100$, the resulting statistics contain $n_j/N$ for
$j = 1,2, \ldots, 120$.

\begin{figure}[!ht]
  \centering
  \subfigure[The empirical PDF versus the PDF of $\mathcal{N}(\frac{1}{2},1/400)$.]{\includegraphics[width=0.45\textwidth]{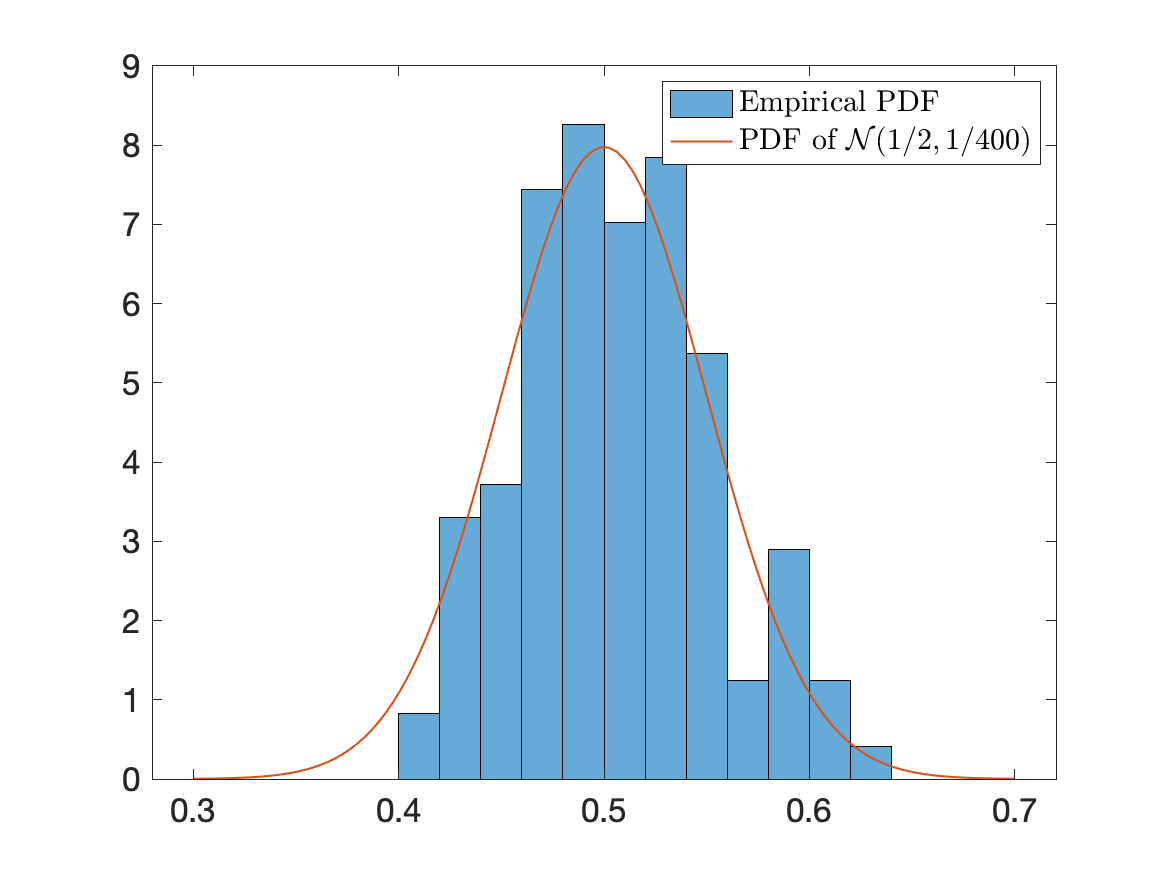}}
  \hspace{0.7em}
  \subfigure[Q–Q plot of $\{n_j/N\}_{j=1}^{120}$
    versus $\mathcal{N}(\frac{1}{2},1/400)$.]{\includegraphics[width=0.45\textwidth]{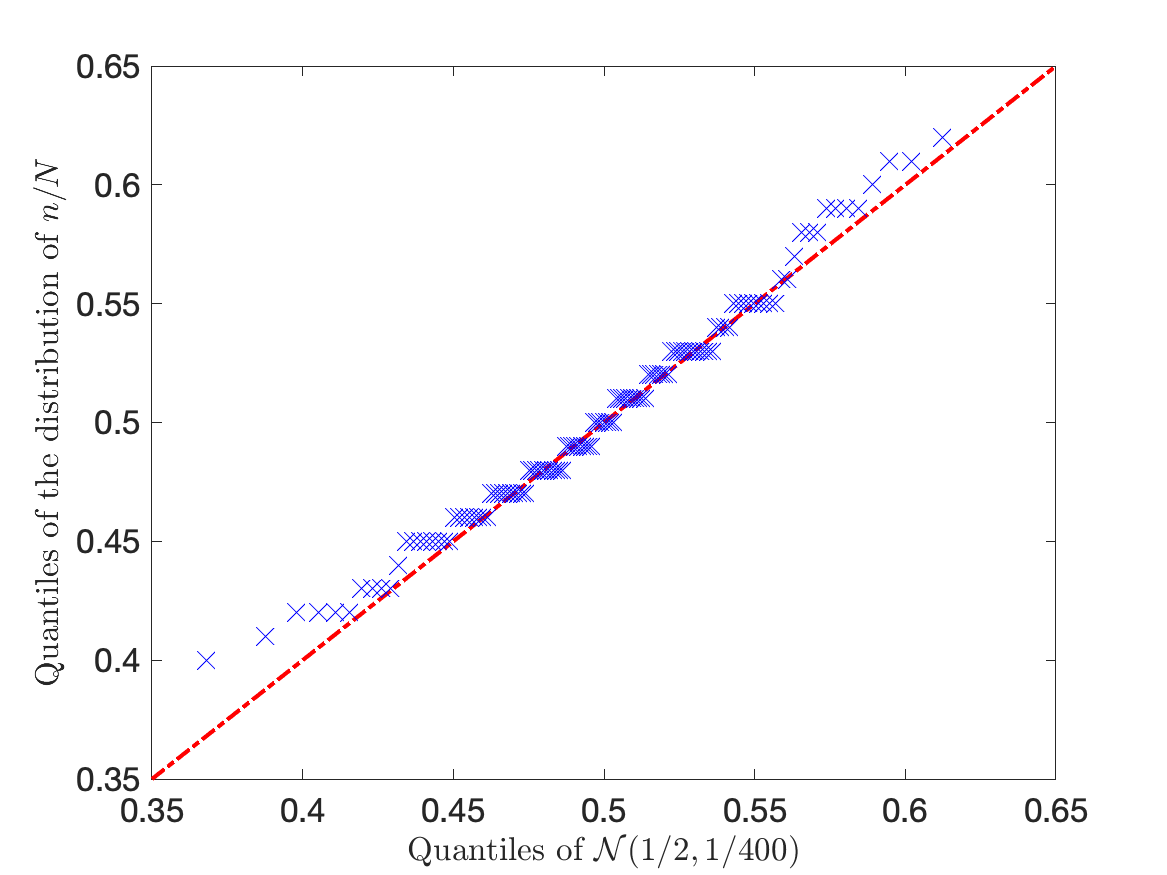} }\label{fig:qq20}
  \caption{Numerical results of the Ginzburg-Landau model when $T=20$. (a): Empirical probability
    density function (PDF) of $\{n_j/N\}_{j=1}^{120}$ compared with the PDF of
    $\mathcal{N}(\frac{1}{2},1/400)$.  (b): Q–Q (quantile-quantile) plot of $\{n_j/N\}_{j=1}^{120}$
    compared with the Q-Q plot of $\mathcal{N}(\frac{1}{2},1/400)$.  }
    \label{fig:isosurface ginzburg 20}
\end{figure}

\begin{figure}[!ht]
    \centering
    \subfigure[The empirical PDF versus the PDF of $\mathcal{N}(\frac{1}{2},1/400)$.]{\includegraphics[width=0.45\textwidth]{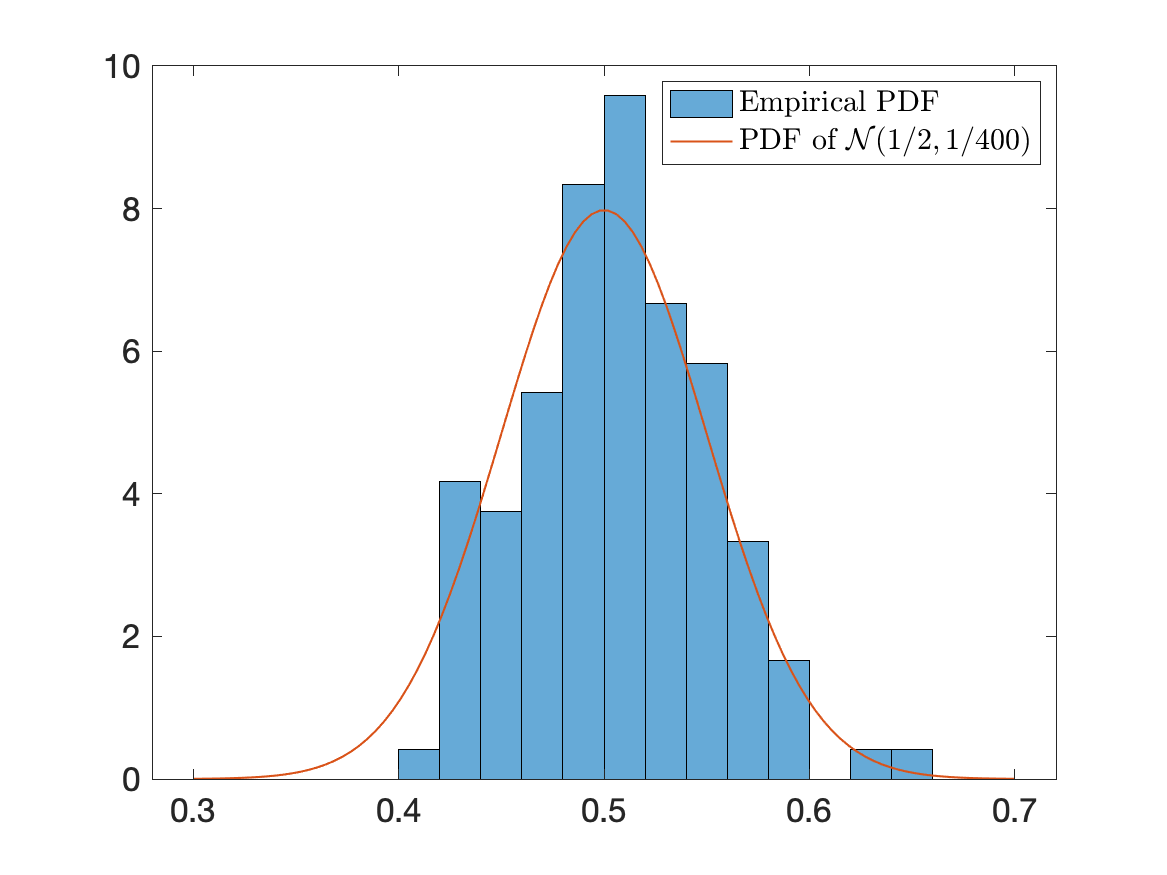}}
    \hspace{0.7em}
    \subfigure[Q–Q plot of $\{n_j/N\}_{j=1}^{120}$
    versus $\mathcal{N}(\frac{1}{2},1/400)$]{\includegraphics[width=0.45\textwidth]{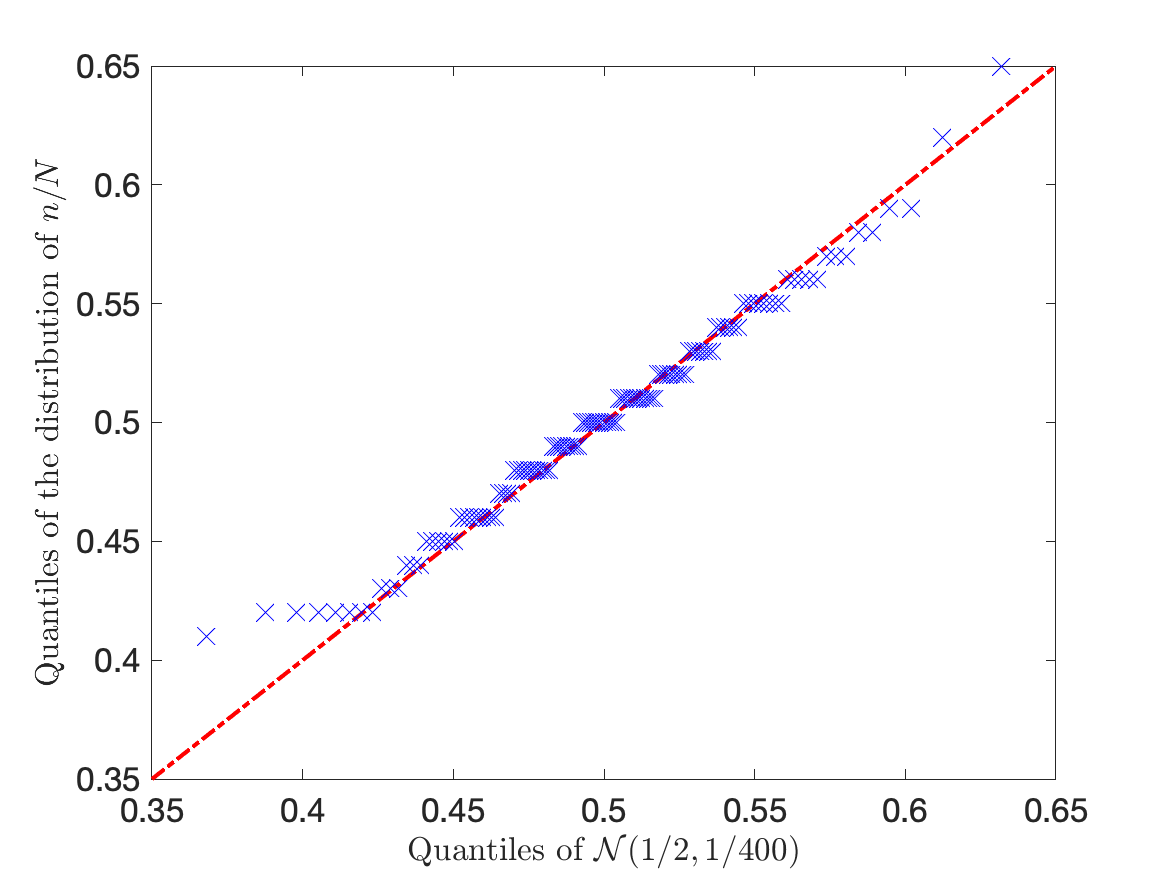} }\label{fig:qq30}
    \caption{Numerical results of the Ginzburg-Landau model when $T=30$. (a): Empirical PDF of $\{n_j/N\}_{j=1}^{120}$ compared with the PDF of $\mathcal{N}(\frac{1}{2},1/400)$. 
    (b): Q–Q (quantile-quantile) plot of $\{n_j/N\}_{j=1}^{120}$ versus $\mathcal{N}(\frac{1}{2},1/400)$.}
    \label{fig:isosurface ginzburg 30}
\end{figure}

The numerical results we get when $T=20$ and $T=30$ are illustrated in Fig.~\ref{fig:isosurface
  ginzburg 20} and Fig.~\ref{fig:isosurface ginzburg 30}, respectively. The histogram of
$\{n_j/N\}_{j=1}^{120}$ is compared with the normal distribution $\mathcal{N}(\frac{1}{2},1/400)$ on
the left, and the Q–Q (quantile-quantile) plot of the distribution of $\{n_j/N\}_{j=1}^{120}$ versus
$\mathcal{N}(\frac{1}{2},1/400)$ is shown on the right. These figures demonstrate that the
distribution of $\{n_j/N\}_{j=1}^{120}$ is in good agreement with the normal distribution
$\mathcal{N}(\frac{1}{2},1/400)$. As mentioned previously in Section~\ref{sec:nonlinear}, we
  can integrate importance sampling \citep{li2019computing,rotskoff2020learning} when dealing with
  metastability, especially when the temperature is relatively low. For instance, we can generate
  the initial state $\x$ in \eqref{eq:nextstep} using importance sampling and then obtain
  $\x_\delta$ via \eqref{eq:nextstep}. 


\section{Conclusion}\label{sec:con}

In this paper, we improve the method in \citep{khoo2019solving} that solves for the neural network
parameterized committor function. In particular, we show that the committor function satisfies an
integral equation based on the semigroup of the Fokker-Planck operator. This integral
  formulation allows us to remove the explicit gradient and handle the boundary conditions
  naturally. The integrals in the variational form of this new equation can be conveniently
approximated via sampling, and the committor function can be solved for using a neural network
parameterization and stochastic gradient descent. The resulting algorithm is shown to be faster
  and less sensitive to the penalty parameter 
  when compared with the approach in
  \citep{khoo2019solving}. The convergence of the training process is guaranteed in the lazy
training regime.

This work suggests a few directions of future research. First on the numerical side, the SDE
\eqref{eq:trajectory} is currently integrated with the Euler-Maruyama method. It will be useful to
explore higher order integration schemes. Our approximation for the stopping time is also rather
primitive and it will be beneficial to explore better decision rules. Second, the proposed method
can be readily applied to other high-dimensional partial differential equations that possess
probabilistic interpretations.


\newpage
\bibliography{ref.bib}

\newpage
\appendix

\section{Proof of Proposition~\ref{prop:semi}}\label{ap:semi}
Recall that the Markov semigroup is defined as follows. 
\begin{defn}\label{def:semi}
The Markov semigroup $\left(P_{t}\right)_{t \geq 0}$ associated with a Markov process $(\x_t)_{t\geq0}$ is defined as
\begin{equation}\label{eq:semidef}
    P_{t} f(\x)=\E^{\x}\left(f\left(\x_{t}\right)\right)=\E\left(f\left(\x_{t}\right) \mid \x_{0}=\x\right), \quad t \geq 0,
    \quad \x \in \Omega,
\end{equation}
where $f:\mathbb{R}^d\rightarrow \mathbb{R}$ is bounded continuous. 
\end{defn}
Recall that $\tau = \tau_{A\cup B}$. While $\tau\wedge \delta$ cannot be directly plugged into Definition~\ref{def:semi} since it is a random variable instead of a constant
, the operator $P_{\tau\wedge \delta}\equiv P$ can still be defined by \eqref{eq:beforesplit}, and by the strong Markov property of the Langevin process $(\x_t)_{t\geq0}$, $(\x_{\tau\wedge \delta})_{\delta\geq0}$ is also a Markov process, thus $(P_{\tau\wedge \delta})_{\delta\geq0}$ is a Markov semigroup.  
Now we proceed to the 
proof with the help of Dynkin's formula described in the following theorem:

\begin{theorem}\citep{dynkin1965markov}

  Consider an It{\^o} diffusion $\{\x_t\}_{t\geq 0}$ defined by the following $d$-dimensional stochastic differential equation (SDE)
\begin{equation}\label{eq:sde}
  d \x_t=  b(\x_t) d t + \sigma(\x_t) d \mathbf{w}_t. 
\end{equation}
Let $\mathcal{A}$ 
be the infinitesimal generator of $\{\x_t\}_{t\geq 0}$, which is defined by 
\begin{equation}\label{eq:infi}
    \mathcal{A} f(\x)=\lim _{t \rightarrow 0+} \frac{P_{t} f(\x)-f(\x)}{t},
\end{equation}
where $f \in D(\mathcal{A}):=\left\{f \in C_{0}\left(\mathbb{R}^{d}\right) \text { such that the limit \eqref{eq:infi} exists}\right\}$, then for a stopping time $\tau$ such that $\E^{\x}\tau < \infty$, 
\begin{equation}\label{eq:dynkin}
    \E^{\x}\left[f\left(\x_{ \tau}\right)\right]=f(\x)+\E^{\x}\left[\int_{0}^{\tau} \mathcal{A} f\left(\x_s\right) \mathrm{d} s\right], 
\end{equation}
for any $f \in C_{0}^{2}\left(\mathbf{R}^{d}\right)$. Moreover, when $\tau$ is the first exit time for a bounded set, \eqref{eq:dynkin} holds for any $f \in C^{2}\left(\mathbf{R}^{d}\right)$
\end{theorem}

For any fixed $\delta>0$, $\tau\wedge \delta$ is a stopping time, and is bounded. Although the committor function $q$ may not be compactly supported, we can still use the formula \eqref{eq:dynkin} with the stopping time $\tau\wedge \delta$ replaced by $\tau\wedge \delta\wedge\tau_{B_r}$, where $B_r$ is the ball with the radius being $r$ and the center being the origin. When $\delta$ is sufficiently small and $r$ sufficiently large, the difference between $P_{\tau\wedge \delta}$ and $P_{\tau\wedge \delta\wedge \tau_{B_r}}$ is negligible in practice, as long as $V$ is a confining potential. 
For this reason, we continue 
to use the notation $P=P_{\tau\wedge \delta}$ in the problem formulation instead of the more cumbersome $P_{\tau\wedge \delta\wedge \tau_{B_r}}$. Therefore 
for the solution $q$ of the equation \eqref{committor}, we have
\begin{equation}\label{eq:replace}
    Pq(\x)=q(\x)+\E^{\x}\left[\int_{0}^{\tau\wedge \delta} -L q\left(\x_s\right) \mathrm{d} s\right] = q(\x), \quad \forall \x \in \Omega\backslash A\cup B, 
\end{equation}
as long as we can prove $\mathcal{A} = -L$. 
In order to show this, we introduce a result in \citep{kallenberg2006foundations}.
\begin{theorem}\citep{kallenberg2006foundations}\label{thm:strong}
Assume that $\x_t$ satisfies the SDE \eqref{eq:sde}. When $a(\x) = \sigma(\x)\sigma(\x)^{\top}$ and $b$ are bounded and Lipschitz continuous, then the generator $\mathcal{A}$ of the semigroup 
associated with \eqref{eq:sde} has the following form. 
\begin{equation}\label{eq:operator}
  \mathcal{A} = \frac{1}{2}\sum_{i,j=1}^da_{ij}(\x)\frac{\partial^2}{\partial x_i\partial
    x_j}+b(\x)\cdot\nabla.
\end{equation}
\end{theorem}

Now we are ready to justify 
\eqref{eq:replace}. 

For the Langevin process \eqref{langevin}, the coefficients are $a_{ij}(\x, t)=\mathbf{1}_{\{i=j\}}1/\beta$, and $b = \grad V$, so the conditions in Theorem \ref{thm:strong} on $a$ are naturally satisfied. Since $\nabla V$ is bounded and Lipschitz on $\mathbb{R}^d$, the conditions in Theorem \ref{thm:strong} on $b$ are satisfied as well. Thus by Theorem \ref{thm:strong}, the generator $\mathcal{A}$ of $(P_t)_{t\geq0}$ is 
\[
\mathcal{A} = \frac{1}{2}\sum_{i,j=1}^d\frac{1}{\beta}\mathbf{1}_{\{i=j\}}\frac{\partial^2}{\partial x_i\partial x_j}+\nabla V(\x)\cdot\nabla = -L.
\]
Plugging this into Dynkin's formula \eqref{eq:dynkin} yields \eqref{eq:justifiedreplace}, which is what we want to prove.

\section{Proof of Proposition~\ref{prop:sym}}\label{ap:sym}

Consider the stationary Langevin process, i.e. $\{\x_t\}_{t\geq0}$ such that  $\x_0\sim\rho$ and $\x_t$ satisfies \eqref{langevin}. It is known that this process is reversible (see for example \citep{pavliotis2014stochastic}). In other words, we have
\[
\E_{\x\sim \rho}u(\x)v(\x_\delta) = \E_{\x\sim \rho}u(\x_\delta)v(\x),
\]
for any $u,v\in L_{\rho}^2(\Omega\backslash A\cup B)$. Now since the law of the trajectory $\{\x_{\tau\wedge t}\}_{t<\delta}$ is the same as $\{\x_t\}_{t<\delta}$ as long as $\delta < \tau$, we have
\[
\E_{\x\sim \rho}u(\x)v(\x_\delta)\mathbf{1}_{\{\delta < \tau\}} = \E_{\x\sim \rho}u(\x_\delta)v(\x)\mathbf{1}_{\{\delta < \tau\}},
\]
which is exactly $\langle u, P^i v\rangle_{\rho} = \langle P^i u, v\rangle_{\rho}$. 

\section{Proof of Proposition~\ref{prop:strcvx}}

  \eqref{eq:pdf} can be rewritten as
\[
P^{i}f(\x) = \int_{\Omega\backslash A\cup B}\!\! \frac{p^i(\x, \mathbf{y})}{\rho(\mathbf{y})}f(\mathbf{y})\rho(\mathbf{y})d\mathbf{y}, 
\]
and when \eqref{eq:HS} holds, we have
\[
\begin{aligned}
&\iint_{\Omega\backslash A\cup B\times\Omega\backslash A\cup B} \left(\frac{p^i(\x, \mathbf{y})}{\rho(\mathbf{y})}\right)^2 \rho(\x)\rho(\mathbf{y})d\x d\mathbf{y}\\
&= \iint_{\Omega\backslash A\cup B\times\Omega\backslash A\cup B}\!\!\!\!\!\!\!\!\!\!\!\!\!\!\!\!\!\!\!\!\!\!\!\!\! p^i(\x, \mathbf{y})\frac{\rho(\x)p^i(\x, \mathbf{y})}{\rho(\mathbf{y})} d\x d\mathbf{y}\\
&=\iint_{\mathrlap{\Omega\backslash A\cup B\times\Omega\backslash A\cup B}} p^i(\x, \mathbf{y})p^i(\mathbf{y}, \x)d\x d\mathbf{y}<\infty, 
\end{aligned}
\]
where in the second equality we used the symmetry of $P^i$. 
Thus $P^i$ is a Hilbert-Schmidt operator on $L_\rho^2(\Omega\backslash A\cup B)$ \citep{dunford1965linear}. Thus it is compact, and every point in its spectrum is an eigenvalue. Since $P^i$ is contractive, all of its eigenvalues are less or equal than $1$, and it suffices to prove that $1$ is not an eigenvalue for $P^i$. If $1$ was an eigenvalue of $P^i$, then there exists a nonzero function $f\in L_\rho^2(\Omega\backslash A\cup B)$ such that $P^i f = f$, and by using \eqref{eq:pdf} $n$ times, we get
\[
\begin{aligned}
f(\x) &= \int_{\mathrlap{\Omega\backslash A\cup B}} p^i(\x, \mathbf{y})f(\mathbf{y})d\mathbf{y} \\
&= \iint_{\mathrlap{\Omega\backslash A\cup B\times \Omega\backslash A\cup B}} p^i(\x, \mathbf{y})p^i(\x, \mathbf{z})f(\mathbf{z})d\mathbf{y}d\mathbf{z}\\
&=\cdots\\
&= \E^{\x}\left(f\left(\x_{n\delta}\right)\mathbf{1}_{\{n\delta<\tau\}}\right), \quad \forall n>0. 
\end{aligned}
\]
But the Langevin process \eqref{langevin} is ergodic, and the equilibrium measure $\rho$ is positive, so $\mathbf{1}_{\{n\delta < \tau\}}\rightarrow0$ almost surely when $n\rightarrow\infty$, resulting in $f = 0$, which contradicts with the assumption that $f$ is nonzero. This contradiction shows that $1$ is not an eigenvalue of $P^i$ and closes the proof.

\end{document}